\documentclass[preprint,12pt]{elsarticle}

\usepackage{amssymb}

\usepackage{amsmath}                    
\usepackage{svg}                        
\usepackage{graphicx}                   
\usepackage{sistyle}                    
\usepackage[defaultlines=3,all]{nowidow}
\usepackage{booktabs}                   
\usepackage{tabularx}                   
\usepackage{adjustbox}                  
\usepackage{textcomp}                   
\usepackage{gensymb}                    
\usepackage{varwidth}                   
\usepackage{supertabular}               
\usepackage{multirow}                   
\usepackage{optidef} 				            
\usepackage{makecell}                   
\usepackage{listings}                   
\usepackage{tikz}                       
\usepackage{pgfplots}                   
\usepackage{mathrsfs}                   
\usepackage{tablefootnote}      
\usepackage{mdframed}
\usepackage{algorithm}
\usepackage{algpseudocode}
\usepackage{amsthm}                     
\usepackage{mathtools}                  

\usepackage{subcaption}
\usepackage{pifont}
\usepackage{soul}

\usepgfplotslibrary{external}
\pgfplotsset{compat=newest}
\usetikzlibrary{positioning,calc,trees,angles,quotes}
\usetikzlibrary{shapes.geometric}
\usetikzlibrary{arrows.meta}
\usetikzlibrary{decorations.pathreplacing}
\usetikzlibrary{intersections}
\usetikzlibrary{fadings}
\usetikzlibrary{calc} 

\definecolor{MATLABblue}{RGB}{0,114,189} 
\definecolor{mygreen}{RGB}{28,172,0}    
\definecolor{mylilas}{RGB}{170,55,241}

\begin{document}

  \begin{frontmatter}


    \title{A variational geometric framework for multi-objective level set topology optimization}

    \author[Oellerich]{Jan Oellerich}
    \author[Yamada]{Takayuki Yamada}

    \affiliation{organization={The University of Tokyo, Graduate School of Engineering, Institute of Engineering Innovation},
              addressline={Yayoi 2-11-16, Bunkyo-ku}, 
              city={Tokyo},
              postcode={113-8656}, 
              country={Japan}}

    \begin{abstract}
    This paper proposes a variational framework for multi-objective level set topology optimization. The approach interprets the level set function as a generalized coordinate of a fictitious material and derives its equation of motion from Hamilton's principle, resulting in a damped wave equation governing the optimization process. The objective functionals are combined using a weighted sum formulation. An analysis of the underlying system structure reveals a geometric interpretation of the problem, shifting the perspective beyond conventional approaches based on purely discrete approximations of the Pareto frontier. Under suitable regularity assumptions, the set of stationary solutions forms a structured subset in objective space, in which the Pareto frontier is locally embedded and the weighting factors act as intrinsic coordinates. This perspective motivates the introduction of a dynamic evolution of the weights, leading to a coupled dynamical system for the level set function and the weighting parameters that enables adaptive exploration of the objective landscape. Numerical results demonstrate that the proposed framework provides a stable and uniform approximation of the Pareto frontier and scales to higher-dimensional objective spaces.
\end{abstract}

    \begin{highlights}
      \item Variational framework for multi-objective level set topology optimization.
      \item Level set evolution derived from Hamilton's principle, resulting in a damped wave equation.
      \item Geometric interpretation of the solution set as a structured subset in objective space.
      \item Dynamic evolution of weighting parameters enabling adaptive exploration of the Pareto frontier.
      \item Improved stability and uniformity compared to conventional weighted sum methods.
    \end{highlights}

  \end{frontmatter}

  \section{Introduction}
\label{sec:Introduction} 
Topology optimization is a powerful method for the design of complex component geometries and offers a higher degree of design freedom compared to sizing \citep{SVANBERG1981,NHA1998} and shape optimization approaches \citep{AZEGAMI2020,WU2023,GUOYONG2025}. Its objective is the identification of an optimal material distribution within a prescribed design domain under given boundary conditions. From a mathematical perspective, this problem is non-trivial and has therefore led to the use of iterative procedures for the approximation of optimal structural layouts.

Among the various methodological classes, boundary based approaches represent a well-established alternative to density based methods, such as SIMP \citep{BENDSOE1999,BENDSOE2004,ZHANG2021}, and evolutionary approaches \citep{XIE1993,HUANG2009,DA2018}, where topological changes are induced through modifications of the material boundary. The level set method is a prominent representative of this class \citep{SETHIAN2000,WANG2003}. Its key concept is the introduction of a scalar function for the implicit representation of the structural boundary, enabling a sharp separation between material and void domains, in contrast to phase-field methods \citep{TAKEZAWA2010,XIE2023,CHEN2026}. The level set method has been successfully applied to a wide range of topology optimization problems \citep{NODA2024,MIYAJIMA2024,HUANG2026}. However, the majority of existing contributions mainly focus on the optimization of a single objective functional. In practical engineering applications, design tasks are frequently governed by multiple, often competing performance criteria. In such settings, the solution is no longer characterized by a single optimal design, but by a set of non-dominated configurations. These configurations constitute the Pareto set, whose image in objective space forms the Pareto frontier, representing trade-offs between conflicting objectives. The accurate approximation of Pareto-optimal solutions in PDE-constrained topology optimization remains challenging and has therefore attracted significant research interest.

A widely used approach to multi-objective topology optimization is the weighted sum (WS) method \citep{ZADEH1963}, in which individual objective functionals are combined into a single scalar functional using non-negative weighting coefficients. While straightforward to implement, this approach restricts the search to convex combinations in objective space and therefore fails to capture non-convex regions of the Pareto frontier. Moreover, the method may lead to a non-uniform distribution of solutions along the Pareto frontier, often resulting in clustering in certain regions. Consequently, various extensions have been proposed to mitigate these concerns.

Izui et al. introduced an aggregative gradient based method in which objective functionals and constraints were evaluated at multiple candidate designs in each iteration \citep{IZUI2015}. The weights were adaptively updated using data envelopment analysis (DEA), thereby assessing the relative performance of the individual candidates. The approach was demonstrated for benchmark problems such as minimum mean compliance and volume minimization within a density based framework. Sato et al. \citep{SATO2016} incorporated a distance-based constraint to prevent clustered point distributions and to enhance solution diversity. A crowded-comparison mechanism was proposed in \citep{SATO2017a} to reduce the number of candidate solutions requiring evaluation. The method was implemented by means of a level set formulation and later applied to the multi-objective design of a no-moving-part valve \citep{SATO2017b} and a microchannel heat sink \citep{SATO2018}. Furthermore, data mining techniques \citep{WU2012,HIWA2014} were employed in \citep{SATO2019} to enable efficient clustering and analysis of structural layouts. An alternative adaptation of the WS method was reported by Chen and Wu who considered mean compliance minimization and maximization of the fundamental eigenvalue \citep{CHEN1998}. Using a density based formulation, the objectives were combined following the compromise programming approach \citep{BORKOWSKI1990}. The method was demonstrated on selected combinations of weighting factors in two numerical examples resulting in a sparse sampling of the objective space. An extension of this method was proposed in \citep{CHEN2000} incorporating fuzzy decision-making techniques to balance the competing influence of compliance and eigenvalue objectives. A similar work was conducted by Luo et al. who employed grey system theory to consider both incomplete and uncertainty information \citep{LUO2006}. Nishiwaki et al. investigated in \citep{NISHIWAKI1998} the multi-objective optimization of compliant mechanisms using the homogenization method \citep{BENDSOE1988,NOGUCHI2021}. The aim was to maximize output flexibility and minimizing mean compliance by considering the ratio of the two objectives. Luo et al. investigated the maximization of mutual energy and minimization of mean compliance in compliant mechanism design using the Rational Approximation of Material Properties (RAMP) method \citep{LUO2005}. The objective functionals were combined similar to the formulation adopted in \citep{CHEN1998}. To ensure numerical stability, additional filtering was employed to suppress checkerboard patterns and one-node connections. Lin et al. presented in \citep{LIN2010} a SIMP-based approach for the design of compliant mechanisms considering minimum mean compliance and maximum mutual energy. The objective functionals were combined via a logarithmic transformation based on physical programming \citep{MESSAC1996,MESSAC2002}. Ryu et. al proposed an adaptive weight determination strategy relying on the construction of a hyperplane in objective space \citep{RYU2019,RYU2021}. The weighting factors were iteratively calculated with respect to the preceding solution. Further work on the incorporation of the WS method can be found in \citep{HOU2025,RYU2026}.

As an alternative scalarization strategy, the Normal Boundary Intersection (NBI) method was developed by Das and Dennis to overcome the limitations of the WS method \cite{DAS1998}. By introducing an additional parametrization of the search direction in objective space, the method enables the generation of evenly distributed candidate solutions along the Pareto frontier. However, Pareto optimality of the resulting solutions is not guaranteed in general. Therefore, Messac et al. proposed a modified formulation incorporating an additional Pareto filtering step to eliminate dominated solutions while maintaining a uniform distribution of sampling points \citep{MESSAC2003,MESSAC2004}. Subsequent refinements were presented in \citep{MARTINEZ2007,MARTINEZ2009}, further improving the distribution properties of the computed Pareto set. Hancock and Mattson proposed the Smart Normal Constraint (SNC) method which employs additional constraints to estimate whether a point is Pareto-efficient \citep{HANCOCK2013}. A modified version of this approach was later presented by Munk et al. \citep{MUNK2018a,MUNK2018b}. The authors coupled the SNC method with the Bi-directional Evolutionary Structural Optimization (BESO) method and demonstrated their approach on minimum mean compliance and on a multi-physics problem involving the bi-objective optimization of a micro fluidic mixer. An application of the SNC method to the optimization of a wing box structural layout was shown by Crescenti et al. using Altair OptiStruct software \citep{CRESCENTI2021}. The method has been later extended for the consideration of more than two objectives \citep{MUNK2019}.

Pareto tracing is another class of methods for computing Pareto-efficient solutions that directly exploits local information along the frontier to generate neighboring solutions and is primarily limited to bi-objective optimization problems. Suresh proposed a Pareto tracing method implemented in MATLAB based on topological differentiation and the SIMP method \citep{SURESH2010}. The approach was demonstrated for a bi-objective problem involving minimum mean compliance and volume minimization. In \citep{TUREVSKY2011}, the method was extended by integrating a fixed-point iteration scheme and applied to the minimum mean compliance problem under multiple load cases. Nathan et al. further developed a Pareto tracing strategy for robust structural optimization \citep{NATHAN2021}, where the expected value and the standard deviation of a structural response were minimized. To improve the diversity and distribution of Pareto-efficient solutions, Polynomial Chaos Expansion (PCE) was employed. Gkaragkounis et al. proposed an adjoint based prediction-correction scheme for Pareto tracing that avoids explicit Hessian evaluations by using Hessian-vector products or Quasi-Newton approximations \citep{GKARAGKOUNIS2021}.

In addition to gradient based and scalarization techniques, population-based evolutionary algorithms have been used in multi-objective optimization. In particular, the Non-dominated Sorting Genetic Algorithm II (NSGA-II) \citep{DEB2002,SCHLEIFER2022} has gained popularity due to its ability to approximate the entire Pareto frontier within a single optimization run without requiring gradient information. However, these approaches typically rely on large populations, which may lead to considerable computational cost. Moreover, the available sensitivity information in structural optimization problems is not explicitly exploited.

A review of the existing literature indicates that multi-objective topology optimization problems are predominantly addressed within density based formulations, whereas comparatively few contributions consider the level set method. Most existing strategies treat multi-objective optimization as an external layer that is superimposed on the underlying topology optimization framework. In the context of level set based approaches, this implies that the intrinsic variational structure governing the evolution of the level set function is not fully accounted for in the multi-objective setting. A consistent extension to multi-objective problems therefore requires a formulation in which the evolution process itself is defined in terms of vector-valued objective functionals. To the best of the authors' knowledge, a systematic formulation of such an extension within the level set framework has not yet been proposed. In \citep{OELLERICH2025}, a variational framework was introduced for the systematic derivation of the evolution equation governing the level set function in the single-objective case. Building upon this foundation, the present work generalizes the variational framework to the multi-objective setting, resulting in a coupled evolutionary system of equations in which both the level set function and the weighting parameters evolve dynamically. Furthermore, the proposed formulation reveals a geometric interpretation of the optimization process in objective space, which naturally gives rise to a simplex-type discretization strategy for the approximation of the Pareto set.

The remainder of this work is organized as follows. Section~\ref{sec:methodology} introduces the fundamental formulation of multi-objective optimization problems in the context of the level set method. On this basis, a coupled evolutionary system is derived and the resulting geometric perspective on the problem is developed. Section~\ref{sec:numericalimplementation} presents the numerical implementation and discretization strategies for the approximation of the Pareto frontier. The influence of numerical parameters and the robustness of the proposed method are investigated in Section~\ref{sec:numericalexperiments} by means of numerical examples. Eventually, Section~\ref{sec:conclusions} summarizes the main findings of this work and outlines potential directions for future research.
  \section{Methodology}
\label{sec:methodology}
\subsection{Level set based topology optimization}
Topology optimization seeks the optimal domain of solid material $\Omega^\prime\subset\mathbb{R}^n,\,n\in\{2,3\}$ that minimizes a given objective functional $J[\mathbf{u}(\Omega)]$. The set of admissible solutions is generally constrained by design restrictions and by the governing equations of the underlying physical system. In this regard, let an isotropic, linear elastic body $\mathcal{B}$ occupy an open domain $\Omega \subseteq \mathrm{D} \subset \mathbb{R}^n$, where $\mathrm{D}$ represents the design domain. The boundary $\partial\Omega$ is assumed smooth and is decomposed into Dirichlet and Neumann parts, $\Gamma^\mathrm{u}$ and $\Gamma^\mathrm{t}$, with the complement $\partial\Omega \setminus (\Gamma^\mathrm{u} \cup \Gamma^\mathrm{t})$. On $\Gamma^\mathrm{t}$ traction forces $t_i$ are imposed while $\Gamma^\mathrm{u}$ is subject to homogeneous Dirichlet boundary conditions. The complement $\partial\Omega \setminus (\Gamma^\mathrm{u} \cup \Gamma^\mathrm{t})$ remains traction-free. Considering indices $i,j,k,l\in\{1,\dots,n\}$ and body forces $b_i$, the equilibrium state in strong form reads
\begin{align}
    \left\{
    \begin{array}{rll}
        -\partial_j(\mathbb{C}_{ijkl}\varepsilon_{kl}(\mathbf{u})) &= b_i &\quad\text{in}\quad\Omega \\[2.5pt]
        u_i &= 0 &\quad\text{on}\quad\Gamma^\mathrm{u} \\[2.5pt]
        (\mathbb{C}_{ijkl}\varepsilon_{kl}(\mathbf{u}))n_j &= t_i &\quad\text{on}\quad\Gamma^\mathrm{t} \\[2.5pt]
        (\mathbb{C}_{ijkl}\varepsilon_{kl}(\mathbf{u}))n_j &= 0 &\quad\text{on}\quad\partial\Omega\setminus(\Gamma^\mathrm{u}\cup\Gamma^\mathrm{t}) 
    \end{array}
    \right.
    \label{eq:equilibrium}
\end{align}
where $\mathbb{C}_{ijkl}$ is the fourth rank elasticity tensor
\begin{align*}
    \mathbb{C}_{ijkl} = \frac{E\nu}{(1+\nu)(1-2\nu)}\delta_{ij}\delta_{kl} + \frac{E}{2(1+\nu)}(\delta_{ik}\delta_{jl} + \delta_{il}\delta_{jk}) 
\end{align*}
incorporating Kronecker delta $\delta$. Parameter $E$ relates to Young's modulus and $\nu$ denotes Poisson's ratio while $\varepsilon_{kl}(\mathbf{u})$ refers to the linearized strain tensor
\begin{align*}
    \varepsilon_{kl}(\mathbf{u}) = \frac{1}{2}(\partial_k u_l + \partial_l u_k)
\end{align*}
being applied to the displacement field $\mathbf{u}$. Considering an additional inequality constraint functional $G[\mathbf{u}]$ and adapting Eq.~(\ref{eq:equilibrium}), the original structural optimization problem is stated as
\begin{customopti}|s|
    {inf}{\Omega\subseteq\mathrm{D}\subset\mathbb{R}^n}{J[\mathbf{u}]=\int_\Omega j(\mathbf{u})\,\mathrm{d}\Omega}{}{}
    \addConstraint{G[\mathbf{u}]}{\leq 0}{}{}
    \addConstraint{-\partial_j(\mathbb{C}_{ijkl}\varepsilon_{kl}(\mathbf{u}))}{=b_i}{\quad\text{in}\quad\Omega}{}
    \addConstraint{u_i}{=0}{\quad\text{on}\quad\Gamma^\mathrm{u}}{}
    \addConstraint{(\mathbb{C}_{ijkl}\varepsilon_{kl}(\mathbf{u}))n_j}{=t_i}{\quad\text{on}\quad\Gamma^\mathrm{t}}{}
    \addConstraint{(\mathbb{C}_{ijkl}\varepsilon_{kl}(\mathbf{u}))n_j}{=0}{\quad\text{on}\quad\partial\Omega\setminus(\Gamma^\mathrm{u}\cup\Gamma^\mathrm{t})}{}.
    \label{eq:problem1}
\end{customopti}
In the next step, Problem~(\ref{eq:problem1}) is transformed into a material distribution problem by employing the characteristic function $\chi_\Omega\in L^\infty(\mathrm{D})$  
\begin{align}
    \chi_\Omega(\mathbf{x}) = 
    \left\{
    \begin{array}{rl}
        1, &\quad\mathbf{x}\in\Omega\cup\partial\Omega\\
        0, &\quad\mathbf{x}\in\mathrm{D}\setminus(\Omega\cup\partial\Omega)
    \end{array}
    \right.,\quad\mathbf{x}\in\mathbb{R}^n
    \label{eq:characteristicfunction}
\end{align}
which allows the original optimization problem to be extended to the design domain $\mathrm{D}$. However, the discontinuous nature of the characteristic function given in Eq.~(\ref{eq:characteristicfunction}) renders the problem ill-posed \citep{ALLAIRE2002}. To overcome this issue, the level set method introduces a higher-dimensional function $\phi:\mathrm{D}\to\mathbb{R},\phi\in H^\kappa(\mathrm{D})$ defined as 
\begin{align}
    \left\{
    \begin{array}{rl}
        1 \geq\phi(\mathbf{x}) > 0, &\quad\mathbf{x}\in\Omega \\
        \phi(\mathbf{x}) = 0, &\quad\mathbf{x}\in\partial\Omega \\
        -1\leq\phi(\mathbf{x})<0, &\quad\mathbf{x}\in\mathrm{D}\setminus(\Omega\cup\partial\Omega)    
    \end{array}
    \right.
    \label{eq:definitionlevelsetfunction}
\end{align}
which enables the implicit description of the structural boundary $\partial\Omega$ by evaluating its isocontour, that is $\phi(\mathbf{x})=0$. Here, $H^\kappa(\mathrm{D})$ denotes the Sobolev space over $\mathrm{D}$ containing weak differentiable functions of order $\kappa\geq 1$. In addition, the characteristic function applied to the level set function in Eq.~(\ref{eq:definitionlevelsetfunction}) can be expressed by a Heaviside function $\Theta(\phi)$:
\begin{align}
    \Theta(\phi)\vcentcolon=
    \left\{
        \begin{array}{rl}
            1,&\quad\text{if}\quad\phi(\mathbf{x}) \geq 0 \\
            0,&\quad\text{otherwise}
        \end{array}
    \right..
    \label{eq:heavisidefunction}
\end{align}
Then, the structural optimization problem for a single objective functional reads
\begin{customopti}|s|
    {inf}{\phi\in H^\kappa(\mathrm{D})}{J[\mathbf{u},\Theta]=\int_\mathrm{D}j(\mathbf{u})\Theta\,\mathrm{d}\Omega}{}{}
    \addConstraint{G[\mathbf{u},\Theta]}{\leq 0}{}{}
    \addConstraint{-\partial_j(\mathbb{C}_{ijkl}\varepsilon_{kl}(\mathbf{u}))}{=b_i}{\quad\text{in}\quad\Omega}{}
    \addConstraint{u_i}{=0}{\quad\text{on}\quad\Gamma^\mathrm{u}}{}
    \addConstraint{(\mathbb{C}_{ijkl}\varepsilon_{kl}(\mathbf{u}))n_j}{=t_i}{\quad\text{on}\quad\Gamma^\mathrm{t}}{}
    \addConstraint{(\mathbb{C}_{ijkl}\varepsilon_{kl}(\mathbf{u}))n_j}{=0}{\quad\text{on}\quad\partial\Omega\setminus(\Gamma^\mathrm{u}\cup\Gamma^\mathrm{t})}{}
    \label{eq:originalproblem}
\end{customopti}
in level set formulation. In contrast to Problem~(\ref{eq:originalproblem}), we are faced with the simultaneous optimization of several competing objective functionals in the case of multi-objective optimization. In our research, we collect these objective functionals in vector form, whereby each objective functional is assigned an individual state variable and thus also an individual equilibrium condition. Consider in total $m\geq 2$ objective functionals indexed by $\alpha\in\{1,\dots,m\}$, then the multi-objective topology optimization problem is formulated as
\begin{customopti}|s|
    {vec\quad inf}{\phi\in H^\kappa(\mathrm{D})}{\mathbf{J}[\mathbf{u}_1,\dots,\mathbf{u}_m,\Theta] = \left(\int_\mathrm{D}j_1(\mathbf{u}_1)\Theta\,\mathrm{d}\Omega, \dots, \int_\mathrm{D}j_{m}(\mathbf{u}_m)\Theta\,\mathrm{d}\Omega\right)^\top}{}{}{}
    \addConstraint{G_\alpha[\mathbf{u}_\alpha,\Theta]}{\leq 0}{}{}
    \addConstraint{-\partial_j(\mathbb{C}_{ijkl}\varepsilon_{kl}(\mathbf{u}_\alpha))}{=b_{\alpha,i}}{\quad\text{in}\quad\Omega}{}
    \addConstraint{u_{\alpha,i}}{=0}{\quad\text{on}\quad\Gamma^\mathrm{u}_\alpha}{}
    \addConstraint{(\mathbb{C}_{ijkl}\varepsilon_{kl}(\mathbf{u}_\alpha))n_j}{=t_{\alpha,i}}{\quad\text{on}\quad\Gamma^\mathrm{t}_\alpha}{}
    \addConstraint{(\mathbb{C}_{ijkl}\varepsilon_{kl}(\mathbf{u}_\alpha))n_j}{=0}{\quad\text{on}\quad\partial\Omega\setminus(\Gamma^\mathrm{u}_\alpha\cup\Gamma^\mathrm{t}_\alpha)}{}.
    \label{eq:originalmultiobjectiveproblem}
\end{customopti}
By means of variational methods, the equilibrium conditions can be expressed in weak form. Therefore, the function spaces
\begin{align*}
    \mathcal{V}_\alpha &= \left\{\mathbf{v}_\alpha=v_{\alpha,i} e_{\alpha,i},v_{\alpha,i}\in H^1(\mathrm{D})\mid\mathbf{v}_\alpha = \mathbf{0}\quad\text{on}\quad\Gamma^\mathrm{u}_\alpha\right\},
\end{align*} 
are first introduced. Then, multiplying by an appropriate test function $\mathbf{v}_\alpha\in\mathcal{V}_\alpha$, integrating over $\mathrm{D}$, and applying the divergence theorem yields the governing equations
\begin{align}
    R_\alpha = \int_\mathrm{D}\mathbb{C}_{ijkl}\varepsilon_{kl}(\mathbf{u}_\alpha)\varepsilon_{ij}(\mathbf{v}_\alpha)\Theta\,\mathrm{d}\Omega - \int_\mathrm{D}b_{\alpha,i}v_{\alpha,i}\Theta\,\mathrm{d}\Omega - \int_{\Gamma^\mathrm{t}_\alpha}t_{\alpha,i} v_{\alpha,i}\,\mathrm{d}\Gamma = 0.
    \label{eq:weakform}
\end{align}
Using Eq.~(\ref{eq:weakform}), Problem~(\ref{eq:originalmultiobjectiveproblem}) can be written as follows:
\begin{customopti}|s|
    {vec\quad inf}{\phi\in H^\kappa(\mathrm{D})}{\mathbf{J}[\mathbf{u}_1,\dots,\mathbf{u}_m,\Theta] = \left(J_1[\mathbf{u}_1,\Theta], \dots, J_m[\mathbf{u}_m,\Theta]\right)^\top}{}{}{}
    \addConstraint{G_\alpha[\mathbf{u}_\alpha,\Theta]}{\leq 0}{}{}
    \addConstraint{R_\alpha[\mathbf{u}_\alpha,\Theta;\mathbf{v}_\alpha]}{= 0}{}{}.
    \label{eq:multiobjectiveproblem}
\end{customopti}
In conventional approaches to single-objective topology optimization, the originally constrained Problem~(\ref{eq:originalproblem}) is typically reformulated as an unconstrained problem by introducing the Lagrangian functional
\begin{align*}
    \mathcal{L} = J[\mathbf{u},\Theta] + \lambda G[\mathbf{u},\Theta] - R[\mathbf{u},\Theta;\mathbf{v}]
\end{align*}
where both the inequality constraint and the equilibrium condition are incorporated into the objective functional, and the inequality constraint is further weighted by the multiplier $\lambda$. A function $\phi^\prime\in H^\kappa(\mathrm{D})$ that satisfies the Karush-Kuhn-Tucker (KKT) conditions
\begin{align}
    \left\{
    \begin{array}{rl}
        \left\langle\partial\mathcal{L}/\partial\phi,\delta\phi\right\rangle &=0 \\
        \lambda G[\mathbf{u},\Theta] &=0 \\
        \lambda &\geq 0 \\
        G[\mathbf{u},\Theta] &\leq 0
    \end{array}
    \right.
    \label{eq:KKTconditions}
\end{align}
can then be considered a solution candidate to the optimization problem. Since a direct determination of $\phi^\prime$ is generally non-trivial, a fictitious time $t$ is introduced, and the level set function is assumed to evolve as $\phi=\phi(\mathbf{x},t)$. The structural boundary is then described by $\phi(\mathbf{x},t)=0$. Differentiation with respect to $t$ leads to the Hamilton-Jacobi equation
\begin{align}
    \partial_t\phi + v_i\partial_i\phi = 0
    \label{eq:hamiltonjacobi}
\end{align}
with $v_i = \mathrm{d}x_i/\mathrm{d}t$, which governs the evolution of the structural boundary. Under appropriate modification, Eq.~(\ref{eq:hamiltonjacobi}) can be reformulated as a reaction-diffusion equation as shown in \cite{YAMADA2010}, allowing the optimal level set function to be approximated iteratively. In the multi-objective context, this approach has been adopted in \citep{SATO2016,SATO2017a,SATO2017b,SATO2019}. 

Note that in contrast to the single-objective case, the multi-objective formulation in Problem~(\ref{eq:multiobjectiveproblem}) does not yield a unique minimizer but rather a set of admissible level set functions. Consequently, the evolution of the level set function can no longer be characterized by a single descent direction. This intrinsic non-uniqueness necessitates the formulation of an equation of motion to describe the dynamics of the level set function in the multi-objective setting. For this purpose, we build upon the variational framework developed in our previous research \cite{OELLERICH2025} in the following.
\subsection{Equation of motion in multi-objective optimization}
\begin{figure}[ht]
    \centering
    \includegraphics[width=0.75\textwidth]{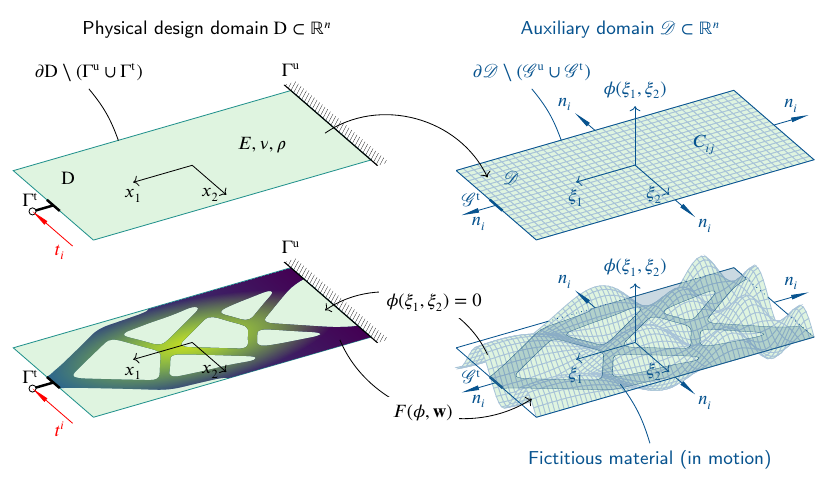}
    \caption{Construction of the auxiliary domain $\mathscr{D}$ containing fictitious material.}
    \label{fig:auxiliarydomain}
\end{figure}

As a starting point, we assume the existence of a solution to Problem~(\ref{eq:multiobjectiveproblem}). From this assumption it follows that the corresponding level set function induces objective values $c_\alpha \in \mathbb{R}$ for each functional $J_\alpha[\mathbf{u}_\alpha,\Theta]$. Consequently, the level set function must also be a solution to the following equations for all $\alpha$:
\begin{align*}
    G_\alpha[\mathbf{u}_\alpha,\Theta] \leq 0,\quad R_\alpha[\mathbf{u}_\alpha,\Theta;\mathbf{v}_\alpha] = 0,\quad J_\alpha[\mathbf{u}_\alpha,\Theta] - c
    _\alpha =0 .
\end{align*}
Let us suppose now that more than one level set function satisfies these conditions. Then all admissible solution candidates can be collected in the set
\begin{align}
    \Phi(\mathbf{c}) = \left\{\bigcap_{\alpha=1}^m\Phi_\alpha^1(\mathbf{c})\right\}\cap\left\{\bigcap_{\alpha=1}^m\Phi_\alpha^2(\mathbf{c})\right\}\cap\left\{\bigcap_{\alpha=1}^m\Phi_\alpha^3(\mathbf{c})\right\}
    \label{eq:solutionset}
\end{align}
with the individual subsets
\begin{align*}
    \Phi^1_\alpha(\mathbf{c}) &\vcentcolon = \{\phi\in H^\kappa(\mathrm{D})\mid G_\alpha[\mathbf{u}_\alpha,\Theta] \leq 0\}, \\
    \Phi^2_\alpha(\mathbf{c}) &\vcentcolon = \{\phi\in H^\kappa(\mathrm{D})\mid R_\alpha[\mathbf{u}_\alpha,\Theta;\mathbf{v}_\alpha] = 0\}, \\
    \Phi^3_\alpha(\mathbf{c}) &\vcentcolon = \{\phi\in H^\kappa(\mathrm{D})\mid J_\alpha[\mathbf{u}_\alpha,\Theta] - c_\alpha = 0\},
\end{align*}
while the $m$ objective functional values are collected in the vector $\mathbf{c}$, that is $\mathbf{c}=(c_1,\dots,c_m)^\top,\mathbf{c}\in\mathbb{J}^m$. 
\begin{figure}[ht]
    \centering
    \includegraphics[width=0.5\textwidth]{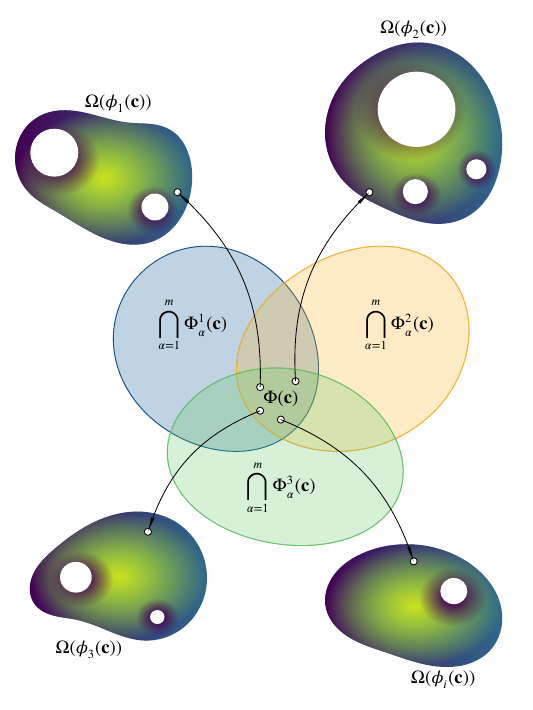}
    \caption{Different topologies $\Omega(\phi_i(\mathbf{c}))$ obtained from $\Phi(\mathbf{c})$.}
    \label{fig:differentsolutions}
\end{figure}
Here, $\mathbb{J}^m\subset\mathbb{R}^m$ denotes the objective functional space defined as
\begin{align*}
    \mathbb{J}^m = \left\{\mathbf{J}=(J_1,\dots,J_m)\in\mathbb{R}^m\mid J_\alpha\in\mathbb{R}\right\}.
\end{align*}
On this basis, we conclude that all admissible level set functions in $\Phi(\mathbf{c})$ also lead to distinct topologies, as qualitatively shown in Fig.~\ref{fig:differentsolutions}. This naturally raises the fundamental question of how to select a representative solution from $\Phi(\mathbf{c})$ in Eq.~(\ref{eq:solutionset}).
\begin{figure}[ht]
    \centering
    \includegraphics[width=0.9\textwidth]{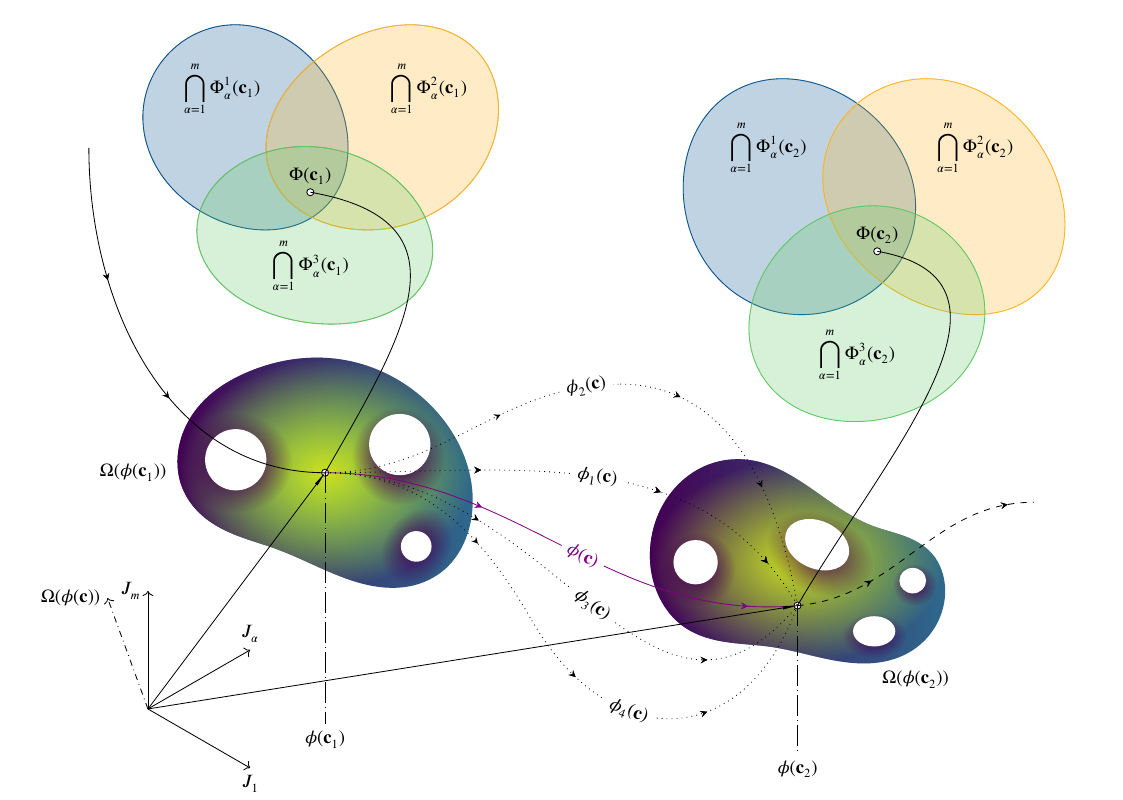}
    \caption{Evolution of the topology $\Omega(\phi(\mathbf{c}))$ in the interval $[\mathbf{c}_1,\mathbf{c}_2]$.}
    \label{fig:evolution}
\end{figure}
However, identifying a unique solution within this set is non-trivial. To address this, we introduce in a first step a neighboring state represented by the vector $\mathbf{c}_2$ and refer to the current state by $\mathbf{c}_1$. With this neighboring state, the according solution set also changes, and with it the admissible level set functions. We note, that attempting to select a specific solution in state $\mathbf{c}_2$ would again face the same ambiguity as in state $\mathbf{c}_1$. For this reason, we adopt a different strategy: we assume the existence of admissible solutions for both states and focus on the transition between them. Although the explicit identification of the solutions is not possible, this assumption allows us to examine the development of an admissible topology from $\mathbf{c}_1$ to $\mathbf{c}_2$. As illustrated in Fig.~\ref{fig:evolution}, multiple trajectories appear conceivable between the two states. To select a meaningful path, we postulate a natural behavior of the system, namely that the evolution of the topology is smooth and free of abrupt changes. In this context, we emphasize once again that the level set function serves as a mathematical representation of the structural boundary. Since only the zero isocontour itself is relevant, it is immaterial in which domain the function is defined. We therefore decouple the level set function from the physical domain $\mathrm{D}$ and consider its evolution in an auxiliary domain $\mathscr{D}$, which is geometrically identical to $\mathrm{D}$ but filled with fictitious material, see Fig.~\ref{fig:auxiliarydomain}. 
This auxiliary domain has no physical reality; it merely provides a conceptual framework for assigning a physically motivated interpretation of the level set function in the following. In particular, we conceptually interpret the fictitious material as a thin elastic film that evolves dynamically. Within this analogy, the level set function can be regarded as a state variable describing the deformation of this film. This perspective enables the definition of energetic measures, such as kinetic and potential energy, and is consistent with the previously postulated assumption of a natural, smooth development of the level set function. To distinguish between the two domains, we furthermore introduce the spatial coordinates $\xi_i$ for the auxiliary domain. Dirichlet and Neumann parts of the boundary $\partial\mathscr{D}$ are denoted by $\mathscr{G}^\mathrm{u}_\alpha$ and $\mathscr{G}^\mathrm{t}_\alpha$. The structural boundary in the physical domain is then represented by the zero level set, $\partial\Omega=\{\xi_i\in\mathbb{R}\mid\phi(\boldsymbol{\xi})=0\}$. In this context, we refer to the previous theoretical observation in Fig.~\ref{fig:evolution} that multiple paths of evolution exist. Due to the assumption of a natural behavior, we postulate that the topology follows a unique path among all admissible trajectories $\phi(\mathbf{c})$ between $\mathbf{c}_1$ and $\mathbf{c}_2$. This finding allows us to invoke the fundamental principle of least action, according to which the associated action $\mathcal{S}$ attains an extremum along the natural evolution of the system. Here, Hamilton's principle represents a particular formulation of the least action principle, stating that the time integral of the Lagrangian $\mathscr{L} = T - U$ attains an extremum, where $T$ denotes the system's kinetic energy and $U$ the potential energy due to the material's internal resistance \citep{BEDFORD2021}. Therefore, Hamilton's principle is applied to our problem setting, while the level set function serves as the generalized coordinate of the system. 
\begin{figure}[ht]
    \centering
    \includegraphics[width=0.9\textwidth]{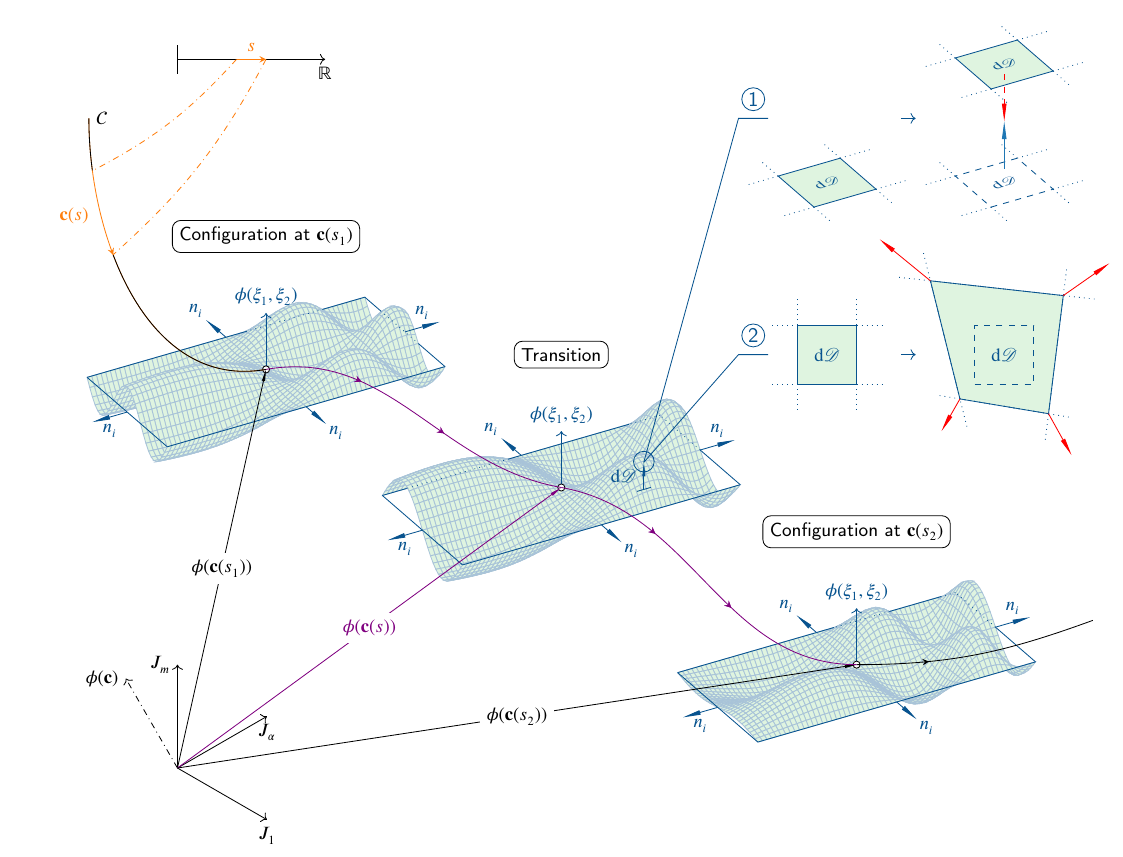}
    \caption{Transition of the level set function from state $\mathbf{c}(s_1)$ to $\mathbf{c}(s_2)$. \ding{192} shows the acceleration of the surface element $\mathrm{d}\mathscr{D}$, \ding{193} shows the compression and tension of a surface element.}
    \label{fig:transition}
\end{figure}
To characterize the evolution more precisely, we describe the development of the topology in terms of a parameterized curve in the objective functional space $\mathbb{J}^m$, where the parameter $s\in\mathbb{R}$ indicates the position on the trajectory, as illustrated in Fig.~\ref{fig:transition}. Let $\mathcal{C}$ denote the smooth curve $\mathcal{C}:[s_1,s_2]\rightarrow\mathbb{J}^m$ with $s\mapsto\mathbf{c}(s)$. We then define an action integral over $\mathcal{C}$ as
\begin{align}
    \mathcal{S}[\phi(\xi_i,\mathbf{c}(s),s)] = \iint_{\mathscr{D}\times\mathcal{C}}\mathscr{L}(\phi,\partial_i\phi,\partial_s\phi,\xi_{i},\mathbf{c}(s),s)\,\mathrm{d}^n \xi\,\mathrm{d}s
    \label{eq:actionintegral}
\end{align}
with operators $\partial_{i} \vcentcolon= \partial/\partial \xi_{i}$ and $\partial_s \vcentcolon= \partial/\partial s$. Now, let $\phi(\mathbf{c}(s))$ denote the true trajectory between the two states $\phi(\mathbf{c}(s_1))$ and $\phi(\mathbf{c}(s_2))$ and let $\delta\phi(\mathbf{c}(s))$ be an admissible variation that vanishes at the boundaries of the interval, i.e. $\delta\phi(\mathbf{c}(s_1))=\delta\phi(\mathbf{c}(s_2))=0$. Furthermore, let $\epsilon$ be a small constant, then the first variation of Eq.~(\ref{eq:actionintegral}) in $\phi$ follows as
\begin{align*}
    \left\langle\frac{\partial\mathcal{S}}{\partial\phi},\delta\phi\right\rangle = \frac{\mathrm{d}}{\mathrm{d}\epsilon}\left[\iint_{\mathscr{D}\times\mathcal{C}}\mathscr{L}(\phi+\epsilon\delta\phi,\partial_{i}(\phi+\epsilon\delta\phi),\partial_s(\phi+\epsilon\delta\phi),\xi_{i},\mathbf{c}(s),s)\,\mathrm{d}^n \xi\,\mathrm{d}s\right]_{\epsilon=0}
\end{align*}
which yields
\begin{align}
    \left\langle\frac{\partial\mathcal{S}}{\partial\phi},\delta\phi\right\rangle = \iint_{\mathscr{D}\times\mathcal{C}}\left\{\frac{\partial\mathscr{L}}{\partial\phi}\delta\phi + \frac{\partial\mathscr{L}}{\partial(\partial_{i}\phi)}\partial_{i}\delta\phi + \frac{\partial\mathscr{L}}{\partial(\partial_s\phi)}\partial_s\delta\phi\right\}\,\mathrm{d}^n \xi\,\mathrm{d}s.
    \label{eq:firstvariation}
\end{align}
The second and third term in Eq.~(\ref{eq:firstvariation}) are treated via integration by parts. Here, we obtain for the second term
\begin{align}
    \iint_{\mathscr{D}\times\mathcal{C}}\frac{\partial\mathscr{L}}{\partial(\partial_{i}\phi)}\partial_{i}\delta\phi\,\mathrm{d}^n \xi\,\mathrm{d}s &= \int_\mathcal{C}\int_{\partial\mathscr{D}}\frac{\partial\mathscr{L}}{\partial(\partial_{i}\phi)}\delta\phi n_{i}\,\mathrm{d}^{n-1} \xi\,\mathrm{d}s - \dots\notag \\
    &\dots - \iint_{\mathscr{D}\times\mathcal{C}}\partial_{i}\left(\frac{\partial\mathscr{L}}{\partial(\partial_{i}\phi)}\right)\delta\phi\,\mathrm{d}^n \xi\,\mathrm{d}s
    \label{eq:intbypartssecond}
\end{align}
and for the third term
\begin{align}
    \iint_{\mathscr{D}\times\mathcal{C}}\frac{\partial\mathscr{L}}{\partial(\partial_s\phi)}\partial_s\delta\phi\,\mathrm{d}^n \xi\,\mathrm{d}s &= \int_\mathscr{D}\left[\frac{\partial\mathscr{L}}{\partial(\partial_s\phi)}\delta\phi\right]_{s=s_1}^{s=s_2}\,\mathrm{d}^n\xi - \dots\notag \\
    &\dots - \iint_{\mathscr{D}\times\mathcal{C}}\partial_s\left(\frac{\partial\mathscr{L}}{\partial(\partial_s\phi)}\right)\delta\phi\,\mathrm{d}^n\xi\,\mathrm{d}s.
    \label{eq:intbypartsthird}
\end{align}
Since the boundary terms vanish in Eq.~(\ref{eq:intbypartssecond}) and Eq.~(\ref{eq:intbypartsthird}) by definition, the first variation reads
\begin{align*}
    \left\langle\frac{\partial\mathcal{S}}{\partial\phi},\delta\phi\right\rangle = \iint_{\mathscr{D}\times\mathcal{C}}\left\{\frac{\partial\mathscr{L}}{\partial\phi} - \partial_{i}\left(\frac{\partial\mathscr{L}}{\partial(\partial_{i}\phi)}\right) - \partial_s\left(\frac{\partial\mathscr{L}}{\partial(\partial_s\phi)}\right)\right\}\delta\phi\,\mathrm{d}^n \xi\,\mathrm{d}s.
\end{align*}
This is satisfied if and only if the term in parentheses vanishes, leading to 
\begin{align}
    \partial_{i}\left(\frac{\partial\mathscr{L}}{\partial(\partial_{i}\phi)}\right) + \partial_s\left(\frac{\partial\mathscr{L}}{\partial(\partial_s\phi)}\right) - \frac{\partial\mathscr{L}}{\partial\phi} = 0.
    \label{eq:eulerlagrangian}
\end{align}
Equation~(\ref{eq:eulerlagrangian}) is commonly known as the Euler-Lagrangian equation and describes in our context the free motion of the fictitious material. However, in the optimization process, the level set function is evolved under the influence of the conditions prevailing in the physical domain. To this end, we incorporate non-conservative terms $Q^\ast$ that steer the evolution of the level set function resulting in the modified Euler-Lagrangian equation
\begin{align}
    \partial_{i}\left(\frac{\partial\mathscr{L}}{\partial(\partial_{i}\phi)}\right) + \partial_s\left(\frac{\partial\mathscr{L}}{\partial(\partial_s\phi)}\right) - \frac{\partial\mathscr{L}}{\partial\phi} = Q^\ast
    \label{eq:eulerlagrangianmodifed}
\end{align}
in strong form. To determine appropriate expressions for the kinetic and potential energy, we theoretically analyze the deformed configuration of the fictitious material, as shown in Fig.~\ref{fig:transition}, by observing the behavior of an infinitesimal surface element $\mathrm{d}\mathscr{D}$. Here, we assume that the fictitious material undergoes acceleration along $\mathcal{C}$, resulting in the kinetic energy 
\begin{align*}
    T(\partial_s\phi) = \frac{1}{2}(\partial_s\phi)^2.
\end{align*}
Regarding the potential energy arising from the internal resistance of the fictitious material, we consider the surface element to experience both compression and tension. Accordingly, we model the potential energy as proportional to the squared gradient of the level set function $\phi$, weighted by a symmetric, positive-definite second-order tensor. To this end, we introduce coefficients $(C_{ij})^2$, which define an effective tensor characterizing the material's resistance to deformation. Thus, the potential energy is defined as
\begin{align*}
    U(\partial_{i}\phi;C_{ij}) = \frac{1}{2}(C_{ij})^2\partial_{i}\phi\partial_{j}\phi.
\end{align*}
The squared coefficients $(C_{ij})^2$ are considered in order to emphasize the analogy of the resulting evolution operator to a generalized wave operator, as will become apparent in the subsequent formulation. Based on these energetic contributions, the Lagrangian density takes the form
\begin{align}
    \mathscr{L} = \frac{1}{2}(\partial_s\phi)^2 - \frac{1}{2}(C_{ij})^2\partial_{i}\phi\partial_{j}\phi.
    \label{eq:lagrangiandensity}
\end{align}
As $\mathscr{L}$ in Eq.~(\ref{eq:lagrangiandensity}) does not explicitly depend on $\phi$, the partial derivatives are computed with respect to the gradient terms. According to Eq.~(\ref{eq:eulerlagrangianmodifed}), we obtain
\begin{align}
    \frac{\partial\mathscr{L}}{\partial(\partial_s\phi)} &= \partial_s\phi, \label{eq:partialderivatives}\\
    \frac{\partial\mathscr{L}}{\partial(\partial_{i}\phi)} &= -\frac{1}{2}\frac{\partial}{\partial(\partial_{i}\phi)}\left((C_{ij})^2\partial_{i}\phi\partial_{j}\phi + (C_{ji})^2\partial_{i}\phi\partial_{j}\phi\right) = -(C_{ij})^2\partial_{j}\phi.\label{eq:partialderivativei}
\end{align}
Note, that these two terms correspond to the contributions governing the free motion of the fictitious material. The non-conservative terms $Q^\ast$ are assumed to originate from a generalized dissipation functional $\mathscr{R}$ via the relation
\begin{align}
    Q^\ast = - \frac{\partial\mathscr{R}}{\partial\phi} - \frac{\partial\mathscr{R}}{\partial(\partial_s\phi)}
    \label{eq:nonconservativeterms}
\end{align}
which is inspired by the structure of a Rayleigh dissipation functional \citep{MOTTAGHI2019}. Here, the first term accounts for changes in the shape of the level set function, whereas the second term captures the rate of change with respect to the curve parameter $s$. In this context, we postulate that the dissipation functional $\mathscr{R}$ can be decomposed into two contributions:
\begin{align*}
    \mathscr{R}[\phi,\partial_s\phi,\xi_i,\mathbf{c}(s),s] = \mathscr{R}_1[\phi,\xi_i,\mathbf{c}(s),s] + \mathscr{R}_2[\partial_s\phi,\xi_i,\mathbf{c}(s),s].
\end{align*}
Since the first term $\mathscr{R}_1$ depends primarily on the current shape of the level set function $\phi$, it is reasonable to quantify this influence on the basis of the given functionals from Problem~(\ref{eq:multiobjectiveproblem}). For this purpose, we propose the following expression for the first part of the dissipation functional:
\begin{align}
    \mathscr{R}_1 &= \sum_{\alpha=1}^{m}\left\{G^\ast_\alpha[\mathbf{u}_\alpha,\Theta;\lambda_\alpha,r_\alpha] + w_\alpha\frac{J_\alpha[\mathbf{u}_\alpha,\Theta] - c_\alpha}{J_\alpha^\ast} - R_\alpha[\mathbf{u}_\alpha,\Theta;\mathbf{v}_\alpha]\right\}.
    \label{eq:dissipationfunctionalfirstterm}
\end{align}
In Eq.~(\ref{eq:dissipationfunctionalfirstterm}), the first summand enforces the inequality constraint, taking into account the individual displacement fields $\mathbf{u}_\alpha$. For numerical implementation, it is denoted by $G^\ast_\alpha$ to indicate that it will be treated using the Augmented Lagrangian method \citep{BURMAN2023} to satisfy the KKT conditions in Eq.~(\ref{eq:KKTconditions}), with each term associated with its own multiplier $\lambda_\alpha$ and penalization factor $r_\alpha$. The objective functionals from Problem~(\ref{eq:multiobjectiveproblem}) are combined into a single aggregated functional based on the WS method, where each objective is scaled by a weighting parameter $w_\alpha\in(0,1)$ to reflect its relative prioritization and additionally related to its initial functional value $J_\alpha^\ast$. This provides a natural scalarization consistent with the variational structure. The weights can be collected to a vector $\mathbf{w}$ belonging to the open simplex defined as
\begin{align}
    \Delta^{m-1} \vcentcolon = \left\{(w_1,\dots,w_m)\in\mathbb{R}^m\;\middle|\;0< w_\alpha < 1,\sum_{\alpha=1}^{m}w_\alpha = 1\right\}.
    \label{eq:opensimplex}
\end{align}
From this, we obtain the partial functional derivative with respect to $\phi$ after separating into implicit and explicit terms \citep{MICHALERIS1994} as follows:
\begin{align}
    \frac{\partial\mathscr{R}}{\partial\phi} &= \sum_{\alpha=1}^{m}\left\{\left(\frac{\partial G^\ast_\alpha}{\partial\mathbf{u}_\alpha} + \frac{w_\alpha}{J_\alpha^\ast}\frac{\partial J_\alpha}{\partial\mathbf{u}_\alpha} - \frac{\partial R_\alpha}{\partial\mathbf{u}_\alpha}\right)\frac{\partial\mathbf{u}_\alpha}{\partial\Theta}\delta(\phi) + \dots \right. \notag\\
    &\left. \dots + \left(\frac{\partial G^\ast_\alpha}{\partial\Theta} + \frac{w_\alpha}{J_\alpha^\ast}\frac{\partial J_\alpha}{\partial\Theta} + \frac{\partial R_\alpha}{\partial\Theta}\right)\delta(\phi)\right\}
    \label{eq:dissipationfunctionalfirsttermderivative}
\end{align}
In this context, $\delta(\phi)$ denotes the Dirac Delta function. From Eq.~(\ref{eq:dissipationfunctionalfirsttermderivative}) we find that the first sum vanishes if each individual summand vanishes. This occurs as soon as the terms in parentheses become zero, which corresponds to the system's underlying adjoint equations. These are given by
\begin{align}
    \mathcal{A}_\alpha:\quad\left\langle\frac{\partial G^\ast_\alpha}{\partial\mathbf{u}_\alpha},\delta\mathbf{u}_\alpha\right\rangle + \left\langle \frac{w_\alpha}{J_\alpha^\ast}\frac{\partial J_\alpha}{\partial\mathbf{u}_\alpha},\delta\mathbf{u}_\alpha\right\rangle = \left\langle\frac{\partial R_\alpha}{\partial\mathbf{u}_\alpha},\delta\mathbf{u}_\alpha\right\rangle,\quad \delta\mathbf{u}_\alpha\in\mathcal{V}_\alpha
    \label{eq:adjointequations}
\end{align}
in weak form. Assuming that these conditions are fulfilled, the partial derivative of the dissipation functional with respect to $\phi$ reduces to
\begin{align}
    \frac{\partial\mathscr{R}}{\partial\phi} = \sum_{\alpha=1}^{m}\left\{\frac{\partial G^\ast_\alpha}{\partial\Theta}  + w_\alpha\left|\frac{1}{J_\alpha^\ast}\frac{\partial J_\alpha}{\partial\Theta}\right| - \frac{\partial R_\alpha}{\partial\Theta}\right\}\delta(\phi) = \sum_{\alpha=1}^m f_\alpha(\phi,w_\alpha)\delta(\phi).
    \label{eq:perturbationterm}
\end{align}
where a single sensitivity contribution is defined as
\begin{align}
    f_\alpha(\phi,w_\alpha) \vcentcolon = \frac{\partial G^\ast_\alpha}{\partial\Theta} + w_\alpha\left|\frac{1}{J_\alpha^\ast}\frac{\partial J_\alpha}{\partial\Theta}\right| - \frac{\partial R_\alpha}{\partial\Theta}.
    \label{eq:perturbationsingle}
\end{align}
Here, the absolute values of the objective functional sensitivities are considered. On this basis, the total sum is conveniently denoted by
\begin{align}
    F(\phi,\mathbf{w}) \vcentcolon= \sum_{\alpha=1}^{m}f_\alpha(\phi,w_\alpha).   
    \label{eq:perturbationsum} 
\end{align}
In the following course of the work, Eq.~(\ref{eq:perturbationterm}) is referred to as the perturbation term containing the individual design sensitivities. In our modeling, the sum of these terms acts as a source of excitation for the fictitious material and thus has a significant influence on the evolution of the level set function, as shown in Fig.~\ref{fig:auxiliarydomain}. In the next step, the second contribution of the dissipation functional is to be specified. Here, we define the expression
\begin{align}
    \mathscr{R}_2[\partial_s\phi,\xi_i,\mathbf{c}(s),s;B] = \frac{1}{2}\iint_{\mathscr{D}\times\mathcal{C}}B(\partial_s\phi)^2\,\mathrm{d}^n \xi\,\mathrm{d}s
    \label{eq:dissipationfunctionalsecondterm}
\end{align}
which corresponds to a kinetic energy dissipation during the evolution multiplied by a damping constant $B > 0$. According to Eq.~(\ref{eq:nonconservativeterms}), differentiation of Eq.~(\ref{eq:dissipationfunctionalsecondterm}) with respect to $\partial_s\phi$ gives
\begin{align}
    \frac{\partial\mathscr{R}}{\partial(\partial_s\phi)} = B\partial_s\phi. 
    \label{eq:derivativedamping}
\end{align}
Finally, insertion of Eqs.~(\ref{eq:partialderivatives}), (\ref{eq:partialderivativei}), (\ref{eq:perturbationsum}), and (\ref{eq:derivativedamping}) into Eq.~(\ref{eq:eulerlagrangianmodifed}) yields
\begin{align*}
    \partial_{ss}\phi - \partial_i((C_{ij})^2\partial_{j}\phi) + B\partial_s\phi = -F(\phi,\mathbf{w})\delta(\phi).
\end{align*}
For compact representations, we introduce the operator 
\begin{align*}
    \square^\prime(\cdot) \equiv \partial_{ss}(\cdot) - \partial_{i}((C_{ij})^2\partial_{j}(\cdot))
\end{align*}
which is a modified version of the D'Alembert operator $\square(\cdot)$ that is commonly used to describe wave characteristics. Following \cite{OELLERICH2025}, we approximate the Heaviside function by means of a hyperbolic tangent function
\begin{align}
    \Theta(\phi;b) \approx \frac{1}{2}(\tanh{(2b \phi)} + 1)
    \label{eq:heavisidefunctionapproximation}
\end{align}
where $b$ controls the width between material and void region. From Eq.~(\ref{eq:heavisidefunctionapproximation}), one obtains through differentiation with respect to $\phi$ the approximated Dirac Delta function
\begin{align}
    \delta(\phi;b) \approx b\,\mathrm{sech}^2(2b\phi).
    \label{eq:diracdeltafunctionapproximation}
\end{align}
A Maclaurin series expansion of Eq.~(\ref{eq:diracdeltafunctionapproximation}) around $\phi=0$, denoting the region between void and material domain, gives
\begin{align*}
    b\,\mathrm{sech}^2(2b\phi) = \sum_{n=1}^{\infty}\left.\frac{1}{n!}\frac{\mathrm{d}^n}{\mathrm{d}\phi^n}\left(b\,\mathrm{sech}^2(2b\phi)\right)\right\vert_{\phi=0}\phi^n.
\end{align*}
By neglecting terms of higher order, the approximated Dirac Delta function reduces to the simple constant $b$, that is
\begin{align*}
    \delta(\phi;b)\approx b.
\end{align*}
Taking into account additional Dirichlet and Neumann boundary conditions, the derived results lead to the following fundamental evolution problem:
\begin{align}
    \left\{
    \begin{array}{rll}
        \square^\prime\phi + B\partial_s\phi &= -F(\phi,\mathbf{w}) b &\quad\text{in}\quad\mathscr{D} \\[9pt]
        \phi &= 1 &\quad\text{on}\quad\mathscr{G}^\mathrm{t}_\alpha \\[9pt]
        \partial_{i}\phi n_i &=0 &\quad\text{on}\quad\partial\mathscr{D}\setminus(\mathscr{G}^\mathrm{u}_\alpha\cup\mathscr{G}^\mathrm{t}_\alpha) \\[9pt]
        \phi\vert_{s=0} &= \phi_0 & \\[9pt] 
        \partial_s\phi\vert_{s=0} &= (\partial_s\phi)_0 &
    \end{array}
    \right.
    \label{eq:evolutionphi}
\end{align}
which fully determines the evolution of the level set function along with the adjoint equations $\mathcal{A}_\alpha$ given in Eq.~(\ref{eq:adjointequations}) and the governing equations $R_\alpha$. It should be emphasized that the proposed formulation does not restrict the evolution to the specific dynamical structure considered here. Alternative evolution laws may be constructed by adopting different energetic contributions or by modifying the balance between conservative and non-conservative contributions as demonstrated in \citep{OELLERICH2025}. The present choice represents one consistent realization within a broader class of admissible dynamical systems.
\subsection{Geometric interpretations}
\begin{figure}[ht]
    \centering
    \includegraphics[width=0.9\textwidth]{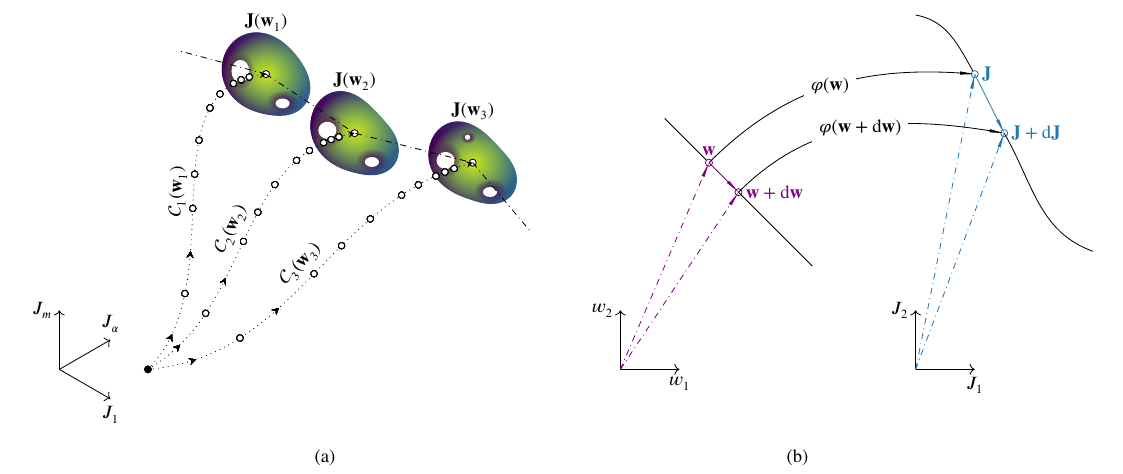}
    \caption{(a) Different solutions obtained for varying weights $\mathbf{w}$, (b) qualitative illustration of the mapping $\varphi$.}
    \label{fig:mapping}
\end{figure}

Equation~(\ref{eq:evolutionphi}) describes a damped wave equation propagating along the curve parameter $s$ in the space of objective functionals $\mathbb{J}^m$. We observe that the motion of a solution candidate is governed by two interacting mechanisms:
\begin{enumerate}
    \item design sensitivities resulting from the optimization problem, which determine a natural descent direction, and
    \item weighting factors $\mathbf{w}$, which act as external control parameters and influence the prioritization of the objective functionals, thereby steering the trajectory in $\mathbb{J}^m$.
\end{enumerate}
In the stationary limit case, the perturbation term $F(\phi,\mathbf{w})$ vanishes and an equilibrium state is reached. Such equilibrium configurations correspond to locally optimal designs under the weighted objective formulation. Different weighting vectors generally lead to different stationary states and thus to different locations of local optima in objective space as qualitatively shown in Fig.~\ref{fig:mapping}a. 
Let $\mathbb{M}\subset\mathbb{J}^m$ denote the set of all stationary objective vectors generated by the evolutionary system. Consequently, this set can be explored systematically by varying the weighting factors. The Pareto frontier forms in this connection a subset of this solution set and is defined with respect to the standard dominance relation in $\mathbb{M}$, i.e.,
\begin{align}
    \mathbb{P} = \left\{\mathbf{J}^\prime\in\mathbb{M} \mid \nexists \mathbf{J}^{\prime\prime} \in \mathbb{M} :\mathbf{J}^{\prime\prime} \preceq \mathbf{J}^\prime \;\wedge\; \mathbf{J}^{\prime\prime} \neq \mathbf{J}^\prime \right\},
    \label{eq:paretoset}
\end{align}
where $\preceq$ denotes the componentwise partial order with respect to minimization. Accordingly, the inclusion $\mathbb{P} \subseteq \mathbb{M} \subset \mathbb{J}^m$ holds. Assuming continuously differentiable objective functionals and local uniqueness of stationary states, small perturbations of the weighting vector $\mathbf{w}$ induce only small changes in the equilibrium configuration and hence small shifts of the associated objective vector in $\mathbb{J}^m$. This motivates the introduction of a continuously differentiable mapping 
\begin{align}
    \varphi : \mathbf{w} \in \Delta^{m-1} \mapsto \mathbf{J} \in \mathbb{J}^m,
    \label{eq:simplexmapping}
\end{align}
which assigns to each admissible weighting vector the corresponding objective vector, see Fig.~\ref{fig:mapping}b. Recall that the open simplex defined in Eq.~(\ref{eq:opensimplex}) is $(m-1)$-dimensional. Under suitable regularity assumptions, in particular assuming sufficiently smooth objective functionals and local uniqueness of stationary states, the mapping $\varphi$ is expected to be locally regular. This suggests that its image is locally at most $(m-1)$-dimensional in $\mathbb{J}^m$, in accordance with the Constant Rank Theorem \citep{TU2007}. In this sense, the solution set $\mathbb{M}$ can locally be regarded as a $(m-1)$-dimensional smoothly parameterized structure. The dependence of stationary solutions on the weighting parameters are then approximated by a first-order Taylor series expansion 
\begin{align*}
    \varphi(\mathbf{w} + \mathrm{d}\mathbf{w})\approx\varphi(\mathbf{w}) + D\varphi(\mathbf{w})\,\mathrm{d}\mathbf{w},     
\end{align*}
so that infinitesimal variations in the weights induce tangential shifts within the solution set. The derivative $D\varphi(\mathbf{w})$ thus characterizes the local sensitivity of stationary states with respect to changes in the prioritization of the objective functionals. This observation gives rise to an interpretation of $\mathbb{M}$ as a differentiable geometric structure embedded in $\mathbb{J}^m$, with the weighting factors $\mathbf{w}$ serving as local coordinates. In this framework, tangential directions and local distance concepts can be defined consistently, allowing the solution set to be interpreted locally as a manifold-like structure. A rigorous manifold-theoretic analysis is reserved for future research. Nevertheless, the present considerations provide the conceptual basis for such an interpretation.
\subsection{Evolution of weightings}
The previous considerations suggest that the stationary solution set $\mathbb{M}$ admits a local geometric structure in objective space, with the weighting vectors $\mathbf{w}$ acting as intrinsic coordinates. During the evolutionary process, however, the topology does not remain fixed but traces a trajectory in $\mathbb{J}^m$. Since the local sensitivity structure and thus the tangent space of $\mathbb{M}$ depend on the current state of the system, a fixed choice of weights imposes a static parametrization on a dynamically evolving geometric object. From this perspective, it appears natural to allow the weighting parameters themselves to evolve adaptively in response to the changing sensitivity landscape. In doing so, the parametrization of the solution set becomes responsive to the local geometric structure, enabling a dynamically consistent exploration of the admissible stationary configurations. Note that the weighting parameters do not evolve freely in $\mathbb{R}^m$, but are geometrically constrained to the open simplex defined in Eq.~(\ref{eq:opensimplex}). For all evolution parameters $s$, they satisfy
\begin{align*}
    \sum_{\alpha=1}^m w_\alpha(s) = 1,\quad 0 < w_\alpha(s) < 1.
\end{align*}
These constraints must be preserved under admissible variations. Rather than enforcing the normalization condition via a Lagrange multiplier in the variational formulation, it is more natural to adopt an intrinsic parameterization of the simplex in which the constraints are satisfied automatically. In such coordinates, the variations become unconstrained, and the geometric structure of $\Delta^{m-1}$ is incorporated directly into the formulation. To this end, we employ the stick-breaking parameterization, originally introduced in the context of Dirichlet processes \citep{SETHURAMAN1994}. Therefore, we introduce the parameter space
\begin{align*}
    \mathbb{Q}^{m-1} = \left\{(q_1,\dots,q_{m-1}) \in \mathbb{R}^{m-1}\;\middle|\;0 < q_\mu < 1 \right\},
\end{align*}
whose elements $\mathbf q \in \mathbb{Q}^{m-1}$ serve as free coordinates. The components $q_\mu$ are indexed by $\mu\in\{1,\dots,m-1\}$. Then, the weights can be defined as smooth functions of $\mathbf{q}$ via
\begin{align}
    \left\{
    \begin{array}{rll}
        w_\nu(\mathbf{q})
        &= (1 - q_\nu)\displaystyle\prod_{\mu=1}^{\nu-1} q_\mu,
        & \quad \nu = 1,\dots,m-1, \\[6pt]
        w_m(\mathbf{q})
        &= \displaystyle\prod_{\mu=1}^{m-1} q_\mu.
    \end{array}
    \right.
    \label{eq:stickbreaking}
\end{align}
This representation is equivalent to the classical stick-breaking construction \cite{ISHWARAN2001,MENA2011,OUYANG2024}. For illustration, in the case of two objective functionals, the weights reduce to
\begin{align*}
    w_1(q_1) = 1 - q_1,\quad w_2(q_1) = q_1,
\end{align*}
while for three objective functionals one obtains
\begin{align*}
    w_1(\mathbf{q}) = 1 - q_1,\quad 
    w_2(\mathbf{q}) = q_1(1 - q_2),\quad 
    w_3(\mathbf{q}) = q_1 q_2.
\end{align*}
The inverse function of Eq.~(\ref{eq:stickbreaking}) accordingly reads
\begin{align}
    q_\mu(\mathbf{w}) = \frac{\sum_{\nu = \mu+1}^{m}w_\nu}{\sum_{\nu=\mu}^{m}w_\nu},\quad \mu=1,\dots,m-1.
    \label{eq:inversestickbreaking}
\end{align}
By construction, positivity and normalization are satisfied identically. The dynamics of the weights can therefore be expressed equivalently in terms of the free coordinates $\mathbf{q}(s)$, which constitute the appropriate variation variables for applying the principle of least action. For $s\in[s_1,s_2],\mathbf{q}\in\mathbb{Q}^{m-1}$, we introduce an action functional
\begin{align}
    \mathcal{S}_q[\mathbf{q}(s)]=\int_{s_1}^{s_2}\mathscr{L}_q(\partial_s\mathbf{q}(s),\mathbf{q}(s),s)\,\mathrm{d}s
    \label{eq:actionfunctionalq}
\end{align}
and calculate the first variation of Eq.~(\ref{eq:actionfunctionalq}) in $\mathbf{q}(s)$ as follows:
\begin{align}
    \left\langle\frac{\partial\mathcal{S}_q}{\partial\mathbf{q}},\delta\mathbf{q}\right\rangle &= \frac{\mathrm{d}}{\mathrm{d}\epsilon}\left[\int_{s_1}^{s_2}\mathscr{L}_q(\partial_{s}(\mathbf{q} + \epsilon\delta\mathbf{q}),\mathbf{q} + \epsilon\delta\mathbf{q})\,\mathrm{d}s\right]_{\epsilon=0} \notag\\
    &= \int_{s_1}^{s_2}\left\{\frac{\partial\mathscr{L}_q}{\partial\mathbf{q}}\cdot\delta\mathbf{q} + \frac{\partial\mathscr{L}_q}{\partial(\partial_s\mathbf{q})}\cdot\partial_s\delta\mathbf{q}\right\}\,\mathrm{d}s.
    \label{eq:firstvariationq}
\end{align}
The second term in Eq.~(\ref{eq:firstvariationq}) is treated via integration by parts which yields
\begin{align*}
    \int_{s_1}^{s_2}\frac{\partial\mathscr{L}_q}{\partial(\partial_s\mathbf{q})}\cdot\partial_s\delta\mathbf{q}\,\mathrm{d}s = \left[\frac{\partial\mathscr{L}_q}{\partial(\partial_s\mathbf{q})}\cdot\delta\mathbf{q}\right]_{s=s_1}^{s=s_2} - \int_{s_1}^{s_2}\partial_s\left(\frac{\partial\mathscr{L}_q}{\partial(\partial_s\mathbf{q})}\right)\cdot\delta\mathbf{q}\,\mathrm{d}s.
\end{align*}
Since the boundary term vanishes, the Euler-Lagrangian equation in strong formulation reads
\begin{align}
    \partial_s\left(\frac{\partial\mathscr{L}_q}{\partial(\partial_s\mathbf{q})}\right) - \frac{\partial\mathscr{L}_q}{\partial\mathbf{q}} = \mathbf{Q}_q
    \label{eq:eulerlagrangianq}
\end{align}
where $\mathbf{Q}_q$ represents non-conservative contributions. To model inertial and restoring effects, we adopt the mechanical ansatz $\mathscr L_q = T_q - U_q$ regarding Eq.~(\ref{eq:eulerlagrangianq}) with
\begin{align*}
    T_q = \frac{1}{2} M_q |\partial_s \mathbf{q}|^2,\quad U_q = \frac{1}{2}K_{q}|\mathbf{q} - \mathbf {q}^\ast|^2 ,
\end{align*}
where $M_q > 0$ and $K_q > 0$ denote effective inertia and stiffness parameters. The potential energy penalizes deviations of $\mathbf{q}$ from a reference state $\mathbf{q}^\ast$ and thereby enables a restoring tendency. The non-conservative term $\mathbf{Q}_q$ is assumed to consist of a damping component to dissipate energy from the system and an information-driven forcing term. In index notation it is proposed as
\begin{align}
    Q_{q,\mu} = -\partial_s\left(\partial_\mu\sum_{\alpha=1}^{m}w_\alpha(\mathbf{q})\frac{J_\alpha}{J_\alpha^\ast}\right) - B_q\partial_s q_\mu
    \label{eq:nonconservativetermsq}
\end{align}
with damping coefficient $B_q > 0$ and operator $\partial_\mu\vcentcolon=\partial/\partial q_\mu$. The first contribution represents the $s$-rate of change of the gradient of the weighted objective functional and does not correspond to classical gradient descent. Instead, it represents the temporal variation of the sensitivity of the weighted objective functional, thereby driving the system according to changes in information rather than absolute objective values. The forcing term vanishes once the sensitivity structure becomes stationary, ensuring asymptotic relaxation of the weight dynamics. In this regard, we extract from Eq.~(\ref{eq:nonconservativetermsq}) the inhomogeneity term
\begin{align*}
    F_{q,\mu}(\phi,\mathbf{q}) = \partial_s\left(\partial_\mu\sum_{\alpha=1}^{m}w_\alpha(\mathbf{q})\frac{J_\alpha}{J_\alpha^\ast}\right),
\end{align*}
which ensures an intrinsic coupling between the parameterized weightings and obtain the resulting $(m-1)$ evolution equations for the weighting parameters:
\begin{align}
    \left\{
        \begin{array}{rl}
            M_q\partial_{ss}q_\mu + B_q\partial_s q_\mu + K_q(q_\mu - q_\mu^\ast) &= -F_{q,\mu}(\phi,\mathbf{q})\quad s\in[s_1,s_2]\\[5pt]
            q_\mu\vert_{s=0} &= (q_{\mu})_0 \\[5pt]
            \partial_s q_\mu\vert_{s=0} &= \partial_s(q_\mu)_0.
        \end{array}
    \right.
    \label{eq:evolutionq}
\end{align}
Note that, Eq.~(\ref{eq:evolutionq}) has the structure of a forced damped oscillator in the parameter space $\mathbb{Q}^{m-1}$. The parameter $M_q$ accounts for inertia, $B_q$ introduces viscous damping, and $K_q$ acts as a stiffness coefficient inducing a restoring tendency toward the reference configuration $\mathbf{q}^\ast$. In the limit $K_q\to\infty$, the conventional WS method is recovered. In this case, the restoring term $K_q(q_\mu - q_\mu^\ast)$ dominates the evolution equation, such that the remaining terms become negligible. The resulting condition $K_q(q_\mu - q_\mu^\ast) \approx 0$ implies $q_\mu = q_\mu^\ast$, i.e., the weights remain fixed at their prescribed reference values. Hence, the classical WS method arises as a limiting case of the proposed weight evolution. The term $F_{q,\mu}(\phi,\mathbf{q})$ represents an external excitation generated by the evolving sensitivity structure of the optimization problem. In this interpretation, the weighting parameters behave as dynamical degrees of freedom whose motion is driven by changes in the objective functionals. In this sense, Eqs.~(\ref{eq:evolutionphi}) and (\ref{eq:evolutionq}) constitute a coupled dynamical system. The two evolution equations are linked through their inhomogeneity terms $F_{q,\mu}(\phi,\mathbf{q})$ and $F(\phi,\mathbf{w}(\mathbf{q}))$, and jointly characterize the trajectory of a solution candidate in objective space.

  \section{Numerical implementation and discretization strategies}
\label{sec:numericalimplementation}
\subsection{Discrete update schemes}
We discretize Eq.~(\ref{eq:evolutionphi}) in the curve parameter $s$ using a finite difference scheme with step size $\Delta s$ and obtain
\begin{align*}
    \frac{\phi^{s+1} - 2\phi^{s} + \phi^{s-1}}{\Delta s^2} + B\frac{\phi^{s+1} - \phi^s}{\Delta s} - \partial_{i}((C_{ij})^2\partial_{j}\phi^{s+1}) = -F^s b.
\end{align*}
Each term $f_\alpha^s$ included in $F^s$ is normalized by a factor $C_\alpha^s$. Rearranging yields
\begin{align}
    (1 + B\Delta s)\phi^{s + 1} - \partial_{i}((C_{ij})^2\partial_{j}\phi^{s + 1})\Delta s^2 = -F^s b\Delta s^2 + (2 + B\Delta s)\phi^s - \phi^{s - 1}
    \label{eq:discretizedequation}
\end{align}
which can be solved for $\phi^{s+1}$ for each iteration $s$. Therefore, let
\begin{align*}
    \Psi = \left\{\delta\phi\in H^\kappa(\mathscr{D})\;\middle|\;\delta\phi=0\quad\text{on}\quad\mathscr{G}_\phi=\bigcup_{\alpha=1}^{m}\mathscr{G}^\mathrm{t}_\alpha\right\},
\end{align*}
then by multiplying Eq.~(\ref{eq:discretizedequation}) with an appropriate test function $\delta\phi\in\Psi$ and integrating over $\mathscr{D}$, the weak form follows as
\begin{align}
    \int_\mathscr{D}(1 + B\Delta s)\phi^{s + 1}\delta\phi\,\mathrm{d}\mathscr{D} - \int_\mathscr{D}\partial_{i}((C_{ij})^2\partial_{j}\phi^{s + 1})\Delta s^2\delta\phi\,\mathrm{d}\mathscr{D} = \dots\notag \\
    \dots = -\int_\mathscr{D}F^s b\Delta s^2\delta\phi\,\mathrm{d}\mathscr{D} + \int_\mathscr{D}\left((2 + B\Delta s)\phi^s - \phi^{s - 1}\right)\delta\phi\,\mathrm{d}\mathscr{D}.
    \label{eq:discreteweakform}
\end{align}
The second term on the left-hand side is treated via integration by parts
\begin{align}
    \int_\mathscr{D}\partial_{i}((C_{ij})^2\partial_{j}\phi^{s + 1})\delta\phi\,\mathrm{d}\mathscr{D} = \int_{\partial\mathscr{D}}\delta\phi (C_{ij})^2\partial_{i}\phi^{s + 1}\,n_{j}\,\mathrm{d}\mathscr{G} - \dots\notag \\
    \dots - \int_\mathscr{D}(C_{ij})^2\partial_{i}\phi^{s + 1}\,\partial_{j}\delta\phi\,\mathrm{d}\mathscr{D}
    \label{eq:discreteintegrationbyparts}
\end{align}
where boundary integrals vanish due to $\delta\phi=0$ on $\mathscr{G}_\phi$ and homogeneous Neumann conditions on the remainder of $\partial\mathscr{D}$. After insertion of Eq.~(\ref{eq:discreteintegrationbyparts}) into Eq.~(\ref{eq:discreteweakform}), the
discretized equation of motion reads
\small{
\begin{align}
    \left\{
        \begin{array}{rl}
            \displaystyle\int_\mathscr{D}(1 + B\Delta s)\phi^{s + 1}\delta\phi\,\mathrm{d}\mathscr{D} + \dots & \\
            \displaystyle \dots + \int_\mathscr{D}(C_{ij})^2\Delta s^2\partial_{i}\phi^{s+1}\partial_{j}\delta\phi \,\mathrm{d}\mathscr{D} &= 
            \displaystyle-\int_\mathscr{D}F^s b\Delta s^2\delta\phi\,\mathrm{d}\mathscr{D} + \dots \\
            &\displaystyle\dots + \int_\mathscr{D}\left((2 + B\Delta s)\phi^s - \phi^{s - 1}\right)\delta\phi\,\mathrm{d}\mathscr{D}\quad\text{in}\quad\mathscr{D} \\
            \phi^{s + 1} &= \displaystyle 1\quad\text{on}\quad\mathscr{G}_\phi=\bigcup_{\alpha=1}^{m}\mathscr{G}_\alpha^\mathrm{t} \\
            \phi\vert_{s=0} &= \phi_0 \\[5pt]
            \partial_s\phi\vert_{s=0} &= \displaystyle\frac{\phi_0 - \phi_{-1}}{\Delta s} 
        \end{array}
    \right.
    \label{eq:discreteevolutionphi}
\end{align}
}\normalsize
which is solved using the finite element method. In this context, let $\mathscr{E}$ be the set of finite elements, then the domain $\mathscr{D}$ is decomposed into the set $\bar{\mathscr{D}}=\bigcup_{e\in\mathscr{E}}V_e$. Furthermore, the continuous fields are approximated by discrete nodal vectors and element-wise matrices
\begin{align*}
    [\mathbf{M}] = \bigcup_{e\in\mathscr{E}}\int_{V_e}N_{\xi} N_{\eta}\,\mathrm{d}V_e,\quad[\mathbf{B}] = \bigcup_{e\in\mathscr{E}}\int_{V_e}\partial_{i} N_{\xi} (C_{ij})^2\partial_{j} N_{\eta}\,\mathrm{d}V_e
\end{align*}
containing shape functions $N_\xi$ and their spatial derivatives $\partial_i N_\xi$ where indices $\xi,\eta$ refer to the total amount of nodes per finite element. Without loss of generality, we set $\Delta s = 1$ between the initial states $s=-1$ and $s=0$, such that only the fields $\{\boldsymbol{\phi}_0\}$ and $\{\boldsymbol{\phi}_{-1}\}$ need to be specified. Finally, we obtain the finite element formulation of Eq.~(\ref{eq:discreteevolutionphi})
\small{
\begin{align}
    \left\{
        \begin{array}{rl}
            \displaystyle\left((1 + B\Delta s)[\mathbf{M}] + \Delta s^2[\mathbf{B}]\right)\{\boldsymbol{\phi}^{s + 1}\} &= \displaystyle[\mathbf{M}](-\{\mathbf{F}^s\} b\Delta s^2 + \dots  \\
            &\dots + (2 + B\Delta s)\{\boldsymbol{\phi}^{s}\} - \{\boldsymbol{\phi}^{s - 1}\})\quad\text{in}\quad\bar{\mathscr{D}} \\
            \{\boldsymbol{\phi}^{s+1}\} &= \displaystyle\{\mathbf{1}\}\quad\text{on}\quad\bar{\mathscr{G}}_\phi = \bigcup_{\alpha=1}^{m}\bar{\mathscr{G}}_\alpha^\mathrm{t} \\[5pt]
            \boldsymbol{\phi}\vert_{s=0} &= \{\boldsymbol{\phi}_0\} \\[9pt]
            \boldsymbol{\phi}\vert_{s=-1} &= \{\boldsymbol{\phi}_{-1}\}
        \end{array}
    \right.
    \label{eq:finiteelementphi}
\end{align}
}\normalsize
with $\bar{\mathscr{G}}_\phi$ as the union set of the discretized boundary sections $\bar{\mathscr{G}}^\mathrm{t}_\alpha$. Next, Eq.~(\ref{eq:evolutionq}) is discretized similarly. Here, we obtain
\begin{align*}
    M_q\frac{q_\mu^{s+1} - 2q_\mu^s + q_\mu^{s-1}}{\Delta s^2} + B_q\frac{q_\mu^{s+1} - q_\mu^s}{\Delta s} + K_q (q_\mu^{s+1} - q_\mu^\ast)= -F_{q,\mu}^s
\end{align*}
which can be explicitly solved for the updated weight component, that is 
\begin{align}
    q_\mu^{s+1} &= \frac{2M_q + B_q\Delta s}{M_q + B_q\Delta s + K_q\Delta s^2}q_\mu^s - \frac{M_q}{M_q + B_q\Delta s + K_q\Delta s^2}q_\mu^{s-1} + \dots\notag \\
    &\dots + \frac{K_q q_{\mu}^\ast\Delta s^2 - F_{q,\mu}^s\Delta s^2}{M_q + B_q\Delta s + K_q\Delta s^2}
    \label{eq:discreteevolutionq}
\end{align}
considering initial conditions $(q_\mu)_{-1}$ and $(q_\mu)_0$. To ensure numerical stability, the weights are constrained to the interval $[\varepsilon_q,1-\varepsilon_q],\varepsilon_q\ll 1$ during the optimization. Concluding, Eqs.~(\ref{eq:finiteelementphi}) and (\ref{eq:discreteevolutionq}) form the discretized coupled evolution problem. The partial differential equations are solved numerically using the software FreeFEM++ \citep{HECHT2012}, combined with the mmg module \citep{DAPOGNY2014,BALARAC2022}, which enables adaptive mesh refinement. The discretization described above provides the numerical framework for solving the coupled evolution equations which govern the motion of one solution candidate. In order to efficiently approximate the Pareto frontier based on the information provided by multiple solution candidates, an adaptive discretization of the weight and objective functional space is proposed in the following.
\subsection{Adaptive Simplex decomposition and algorithmic implementation}
\begin{figure}[ht]
    \centering
    \includegraphics[width=0.9\textwidth]{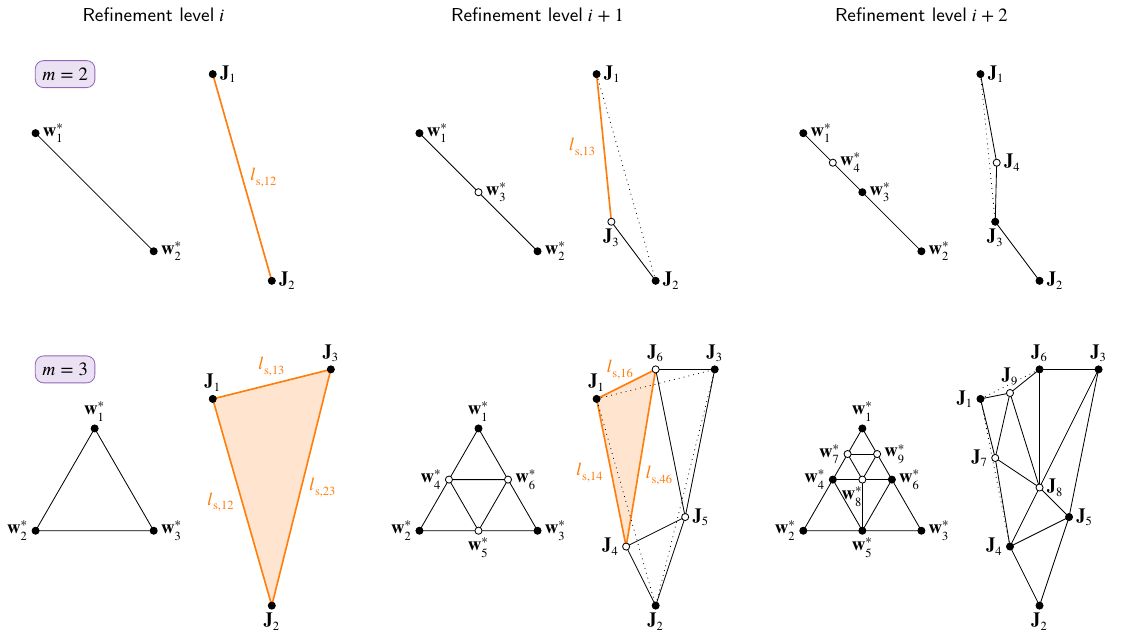}
    \caption{Qualitative illustration of the ASD method; poor simplices marked in orange.}
    \label{fig:simplexdecomposition}
\end{figure}

According to the previously introduced geometric interpretation, it is natural to discretize the weight space $\Delta^{m-1}$ by means of a simplicial decomposition. This induces a corresponding discretization of the mapping $\varphi$ given in Eq.~(\ref{eq:simplexmapping}) from the weight space to the objective functional space. Each simplex in the discretized weight space is therefore mapped to a simplex in the objective space spanned by the respective objective vectors. In order to measure the quality of the approximation of the Pareto frontier, the simplices in the objective space are evaluated with respect to their edge lengths $l_\mathrm{s}$. For this purpose, the normalized objective values are considered to account for differences in scale between the objectives. If the maximum edge length of a simplex exceeds a prescribed threshold $l_\mathrm{s,max}$, the simplex is marked as poor and scheduled for refinement in the subsequent iteration. Under the assumption of a continuous mapping $\varphi$ between $\Delta^{m-1}$ and $\mathbb{J}^m$, weights located in the interior of a simplex in weight space correspond to solutions located within the associated simplex in the objective space. Based on this observation, new weights are generated using the barycenters of the edges of the marked simplices. In this way, the approximation of the solution set is successively refined. To quantify the global quality of the approximation, the mean edge length $\bar{\mu}(l_\mathrm{s})$ of all simplices in the objective space is taken into account. The refinement procedure is terminated once this quantity falls below the prescribed tolerance $l_\mathrm{s,max}$ or after reaching a certain amount of refinement steps. In the following course of the work this approach is referred to as the Adaptive Simplex Decomposition (ASD) method. A qualitative illustration of the refinement strategy is shown in Fig.~\ref{fig:simplexdecomposition}. 
For the case of more than two objective functionals, the simplicial decomposition is constructed using a Delaunay triangulation. This preserves topological neighborhood relationships and allows the approach to be extended, in principle, to higher-dimensional objective spaces \citep{LAWSON1986}, thereby supporting the scalability of the proposed method. We note that a related refinement strategy was proposed by Kim and de Weck \citep{KIM2005,KIM2006}. In the adaptive WS method, additional constraints in the objective space were introduced to guide the search between adjacent Pareto points. While this strategy is effective for general optimization frameworks, the incorporation of such constraints is non-trivial in a level set based formulation, since it would require the introduction of additional Lagrange multipliers and modify the underlying evolution equation. In contrast, the approach proposed in this work does not alter the original optimization problem. The adaptive sampling of the Pareto frontier is achieved solely through a refinement of the discretization of the weight space combined with an adaptive evolution of the weights, which makes the method naturally compatible with the level set framework. The overall structure of the proposed algorithm is summarized in Algorithm~\ref{alg:algorithm}. 
\begin{algorithm}[ht]
\caption{Adaptive simplex decomposition method}
\label{alg:}
\begin{algorithmic}[1] 
\State Define input settings
\State Initialize the discrete reference weight set $\mathcal{W}_h^\ast = \{\mathbf{w}_i^\ast\}$
\While{mean simplex edge length $\bar{\mu}(l_\mathrm{s}) > l_\mathrm{s,max}$}
    \For{each reference weight $\mathbf{w}_i^\ast \in \mathcal{W}_h^\ast$}
        \State Initialize level set functions $\phi_0,\phi_{-1}$
        \For{$s = 1,\dots,N_{\max}$}
            \State Update weightings $\mathbf{q}^s$ and convert to $\mathbf{w}^s$
            \State Solve state problems $R_\alpha$ and adjoint problems $\mathcal{A}_\alpha$
            \State Evaluate objectives $J_\alpha$ and constraints $G_\alpha$
            \State Compute the aggregated perturbation $F(\phi^s,\mathbf{w}^s)$
            \State Update level set function $\phi^{s+1}$
            \If{stationarity criterion is satisfied}
                \State \textit{break}
            \EndIf
        \EndFor
        \State Store solution candidates
    \EndFor
    \State Update solution register
    \State Perform simplex construction (\textsc{Delaunay} for $m>2$)
    \State Refine poor simplexes and update reference weight set $\mathcal{W}_h$
\EndWhile
\State Remove closely spaced points and filter for Pareto efficient solutions
\end{algorithmic}
\label{alg:algorithm}
\end{algorithm} The procedure can be divided into an inner and an outer loop. After defining the input parameters and initializing the discretized reference weight set $\mathcal{W}_h^\ast$, the outer loop performs an adaptive refinement of the weight space until the termination criterion based on the mean simplex edge length is satisfied. The inner loop comprises the numerical solution of the coupled evolution problem for the level set function and the weights. For each reference weight, the state and adjoint problems are solved and the aggregated perturbation field is evaluated in order to update the level set function. The weight parameters are updated simultaneously according to the discrete evolution scheme. Once the stationarity criterion is satisfied, the corresponding solution candidate is stored in a global solution register. Based on the collected solutions, a simplicial decomposition of the objective space is constructed. Poor simplices are identified according to the edge-length criterion, and additional reference weights are introduced to refine the discretization of the weight space. Filtering with respect to Pareto efficiency is performed in the final step. This allows the non-dominated subset according to Eq.~(\ref{eq:paretoset}) to be extracted from the adaptively reconstructed geometric structure of the solution set in the objective space, which acts as the underlying carrier of the Pareto frontier. Approximating this geometric structure rather than the Pareto frontier alone provides additional insights into the optimization problem, for instance regarding the local curvature of the mapping between weight space and objective space.

  \section{Numerical experiments}
\label{sec:numericalexperiments}
In this section, the proposed framework is demonstrated using several numerical examples. To avoid an excessive number of indices, we adopt the convention that quantities and their indices are assumed to correspond to the respective example under consideration. Furthermore, index $\alpha$ is only included if a quantity differs between the single optimization problems. If a parameter is identical for all problems, the index is omitted.
\subsection{Bi-objective mean compliance minimization}
\begin{figure}[ht]
    \centering
    \includegraphics[width=0.9\textwidth]{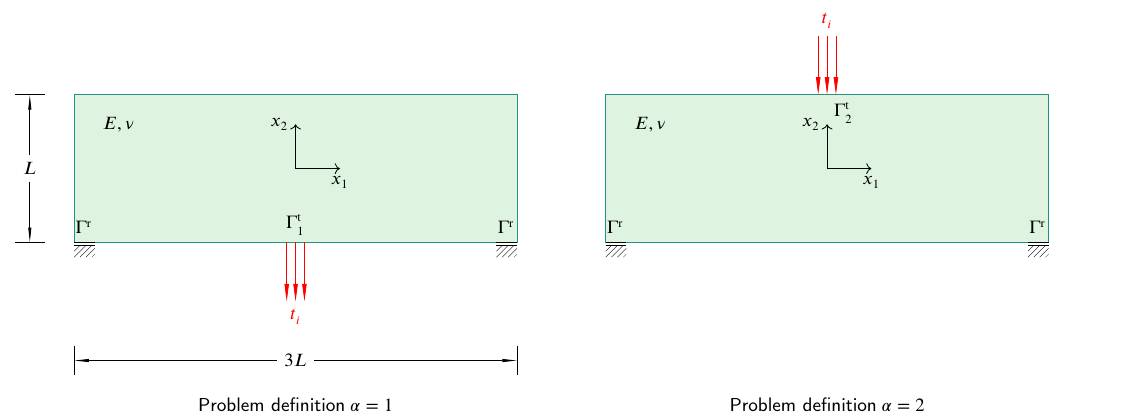}
    \caption{Bi-objective minimum mean compliance problem: simply supported girder.}
    \label{fig:problemGRD}
\end{figure}

We study a minimum mean compliance problem in a girder structure ($\alpha\in\{1,2\}$), where the material domain $\Omega$ is constrained in normal direction at $\Gamma^\mathrm{r}$ as shown in Fig.~\ref{fig:problemGRD}. 
Considering traction forces $t_{i}$ being applied at $\Gamma^\mathrm{t}_\alpha$ and neglecting body forces, the optimization problem is formulated as  
\begin{customopti}|s|
    {vec\quad inf}{\phi\in H^2(\mathscr{D})}{
    \mathbf{J}[\mathbf{u}_1,\mathbf{u}_2] = \left(J_1[\mathbf{u}_1],J_2[\mathbf{u}_2]\right)^\top}{}{}{}
    \addConstraint{G[\Theta] = \frac{1}{V_0}\int_\mathrm{D}\Theta\,\mathrm{d}\Omega - V_\mathrm{f}}{\leq 0}{}{}
    \addConstraint{-\partial_j(\mathbb{C}_{ijkl}\varepsilon_{kl}(\mathbf{u}_\alpha))}{=0}{\quad\text{in}\quad\Omega}{}
    \addConstraint{u_{\alpha,i} n_i}{=0}{\quad\text{on}\quad\Gamma^\mathrm{r}}{}
    \addConstraint{(\mathbb{C}_{ijkl}\varepsilon_{kl}(\mathbf{u}_\alpha))n_j}{=t_i}{\quad\text{on}\quad\Gamma^\mathrm{t}_\alpha}{}
    \addConstraint{(\mathbb{C}^{ijkl}\varepsilon_{kl}(\mathbf{u}_\alpha))n_j}{=0}{\quad\text{on}\quad\partial\Omega\setminus(\Gamma^\mathrm{r}\cup\Gamma^\mathrm{t}_\alpha)}{}
    \label{eq:biminimummeancompliance}
\end{customopti}  
with objective functionals  
\begin{align*}
    J_\alpha[\mathbf{u}_\alpha] = \int_{\Gamma^\mathrm{t}_\alpha} t_i u_{\alpha,i} \,\mathrm{d}\Gamma.
\end{align*}
Here, parameter $V_0$ denotes the volume of the design domain. From the boundary conditions of Problem~(\ref{eq:biminimummeancompliance}), the governing equations are derived in weak form as  
\begin{align}
    R_\alpha:\quad \int_\Omega \mathbb{C}_{ijkl}\varepsilon_{kl}(\mathbf{u}_\alpha) \varepsilon_{ij}(\mathbf{v}_\alpha) \,\mathrm{d}\Omega - \int_{\Gamma^\mathrm{t}_\alpha} t_i v_{\alpha,i} \,\mathrm{d}\Gamma = 0, \quad \mathbf{v}_\alpha \in \mathcal{V}_\alpha.
    \label{eq:governingequationomega}
\end{align}  
To extend Eq.~(\ref{eq:governingequationomega}) to the full design domain $\mathrm{D}$, we adopt an Ersatz material approach \citep{CHOI2011}. Therefore, an auxiliary function  
\begin{align}
    \tau(\Theta) = (1-d)\Theta^a + d, \quad d \ll 1, \, a > 1,
    \label{eq:tau}
\end{align}  
is introduced to express the elasticity tensor as $\mathbb{C}_{ijkl,\phi}(\Theta) = \tau(\Theta)\mathbb{C}_{ijkl}$. Noting that $\Theta \in \{0,1\}$ according to Eq.~(\ref{eq:heavisidefunction}), the identities $\Theta^a(\phi) \equiv \Theta^{a-1}(\phi) \equiv \Theta(\phi)$ hold, which simplifies the derivative $\partial\tau/\partial\Theta = a(1-d)\Theta^{a-1}$ of Eq.~(\ref{eq:tau}) with respect to $\Theta$ to  
\begin{align}
    \frac{\partial \tau}{\partial \Theta} = a(1-d)\Theta.
    \label{eq:partialtau}
\end{align}  
Thus, the governing equations in weak form over the entire design domain are written as  
\begin{align}
    R_\alpha:\quad \int_\mathrm{D} \tau(\Theta) \mathbb{C}_{ijkl} \varepsilon_{kl}(\mathbf{u}_\alpha) \varepsilon_{ij}(\mathbf{v}_\alpha) \,\mathrm{d}\Omega - \int_{\Gamma^\mathrm{t}_\alpha} t_i v_{\alpha,i} \,\mathrm{d}\Gamma = 0,\quad\mathbf{v}_\alpha\in\mathcal{V}_\alpha
    \label{eq:GRDgoverningequations}
\end{align}  
in level set formulation. Since the constraint functional does not depend on $\mathbf{u}_\alpha$, the respective adjoint equations follow directly from the remaining terms:  
\begin{align}
    \mathcal{A}_\alpha:\quad \int_{\Gamma^\mathrm{t}_\alpha}\frac{w_\alpha}{J_\alpha^\ast} t_i \delta u_{\alpha,i} \,\mathrm{d}\Gamma = \int_\mathrm{D} \tau(\Theta) \mathbb{C}_{ijkl}\varepsilon_{kl}(\delta \mathbf{u}_\alpha) \varepsilon_{ij}(\mathbf{v}_\alpha) \,\mathrm{d}\Omega,\quad\delta\mathbf{u}_\alpha\in\mathcal{V}_\alpha.
    \label{eq:GRDadjointequations}
\end{align}  
The partial functional derivative of the inequality constraint functional with respect to $\Theta$ reads
\begin{align*}
    \left\langle\frac{\partial G}{\partial\Theta},\delta\Theta\right\rangle &= \frac{\mathrm{d}}{\mathrm{d}\epsilon}\left[\frac{1}{V_0}\int_\mathrm{D}(\Theta + \epsilon\delta\Theta)\,\mathrm{d}\Omega\right]_{\epsilon=0} = \frac{1}{V_0}\int_\mathrm{D}\delta\Theta\,\mathrm{d}\Omega.
\end{align*}
From this, we identify the partial derivative 
\begin{align}
    \frac{\partial G}{\partial\Theta} = \frac{1}{V_0}.
    \label{eq:GRDpartialderivativeconstraint}
\end{align}
Similarly, the partial functional derivatives of the governing equations follow as
\begin{align*}
    \left\langle\frac{\partial R_\alpha}{\partial\Theta},\delta\Theta\right\rangle = \frac{\mathrm{d}}{\mathrm{d}\epsilon}\left[\int_\mathrm{D}\tau(\Theta + \epsilon\delta\Theta)\mathbb{C}_{ijkl}\varepsilon_{kl}(\mathbf{u}_\alpha)\varepsilon_{ij}(\mathbf{v}_\alpha)\,\mathrm{d}\Omega - \int_{\Gamma^\mathrm{t}_\alpha}t_i v_{\alpha,i}\,\mathrm{d}\Gamma\right]_{\epsilon=0}
\end{align*}
which yields
\begin{align}
    \frac{\partial R_\alpha}{\partial\Theta} = \frac{\partial\tau}{\partial\Theta}\mathbb{C}_{ijkl}\varepsilon_{kl}(\mathbf{u}_\alpha)\varepsilon_{ij}(\mathbf{v}_\alpha).
    \label{eq:GRDpartialderivativegoverning}
\end{align}
Using Eqs.~(\ref{eq:partialtau}), (\ref{eq:GRDpartialderivativeconstraint}) and (\ref{eq:GRDpartialderivativegoverning}), we subsequently obtain the respective perturbation terms based on Eq.~(\ref{eq:perturbationsingle}) which are required to update the level set function:
\begin{align}
    f_1 &= \frac{\lambda}{2V_0} - \frac{1}{C_1^s}\frac{\partial\tau}{\partial\Theta}\mathbb{C}_{ijkl}\varepsilon_{kl}(\mathbf{u}_1)\varepsilon_{ij}(\mathbf{v}_1), \label{eq:GRDperturbationterm1}\\
    f_2 &= \frac{\lambda}{2V_0} - \frac{1}{C_2^s}\frac{\partial\tau}{\partial\Theta}\mathbb{C}_{ijkl}\varepsilon_{kl}(\mathbf{u}_2)\varepsilon_{ij}(\mathbf{v}_2).
    \label{eq:GRDperturbationterm2}
\end{align}
The inequality constraint is active in both optimization problems, therefore a single Lagrangian multiplier $\lambda$ is introduced, which is equally distributed between the two problems and appears as $\lambda/2$ in each case. Due to the self-adjoint nature of the problem, the adjoint variable is identical to the corresponding state variable. Therefore, the sum of the perturbation terms given by Eq.~(\ref{eq:GRDperturbationterm1}) and (\ref{eq:GRDperturbationterm2}) can be reformulated according to Eq.~(\ref{eq:perturbationsum}) as:
\begin{align}
    F = \frac{\lambda}{V_0} - \sum_{\alpha=1}^{2}\left\{\frac{1}{C_\alpha^s}\frac{\partial\tau}{\partial\Theta} \frac{w_\alpha}{J_\alpha^\ast}\mathbb{C}_{ijkl}\varepsilon_{kl}(\mathbf{u}_\alpha)\varepsilon_{ij}(\mathbf{u}_\alpha)\right\}. 
    \label{eq:GRDperturbationterm}
\end{align}
Furthermore, the normalization factor
\begin{align}
    C_\alpha^s = \frac{1}{w_\alpha V_0}\int_\mathrm{D}\left\Vert\frac{\partial R_\alpha}{\partial\Theta}\right\Vert_{L^2(\mathrm{D})}\,\mathrm{d}\mathrm{D}
    \label{eq:GRDnormalizationfactor}
\end{align}
is incorporated to balance the magnitude of the perturbation terms associated with the different objectives. The evolution problem is now fully defined by Eqs.~(\ref{eq:GRDgoverningequations}), (\ref{eq:GRDadjointequations}), and (\ref{eq:GRDperturbationterm}) and forms the basis for investigating the impacts of the numerical parameters.
\subsubsection{Influence of the numerical parameters}
To study the influences of the individual parameters, a reference configuration is defined using the parameter combination $M_q=B_q=K_q=1$ and initial conditions defined by the ratio $\bar{q}_1=(q_1)_0/q_1^\ast=0$. If a parameter is varied in the numerical studies, all remaining parameters are kept fixed. In all experiments, Young's modulus is set to $E=1\,\mathrm{N/mm^2}$ and Poisson's ratio to $\nu=0.3$. The volume fraction is set to $V_\mathrm{f}=0.45$. A traction load of magnitude $|\mathbf{t}|=1\,\mathrm{N}$ is applied, and the length of the domain is given as $L=1\,\mathrm{mm}$. The parameters of the Ersatz material approach are chosen as $a=3$ and $d=1\times 10^{-3}$. The initial layout is a solid domain, i.e. $\phi_{-1}=\phi_0=1$. Furthermore, the evolution parameters are set to $C_{11}=C_{22}=0.014$ and $b=1$ with a damping coefficient of $B=0.001$. The finite element size varies between $h_{\min}=0.5\times10^{-4}$ and $h_{\max}=1\times10^{-2}$. To compare the adaptive evolution of the weights with the classical WS method, we examine the solution trajectory in the objective space for the reference weights $\mathcal{W}_h^\ast=\{(0.9,0.1)\}$ which correspond to $q_1^\ast = 0.9$ following Eqs.~(\ref{eq:stickbreaking}) and (\ref{eq:inversestickbreaking}). 

\begin{figure}[ht]
    \centering
    \includegraphics[width=0.9\textwidth]{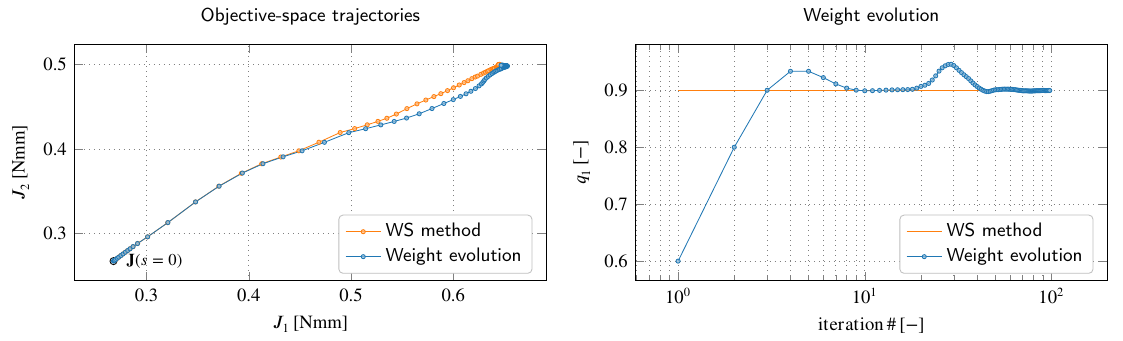}
    \caption{Comparison of the WS method and the proposed weight evolution scheme in terms of the objective-space trajectories and the evolution of the weight parameter $q_1$.}
    \label{fig:comparisonmethod}
\end{figure}

Figure~\ref{fig:comparisonmethod} shows the resulting trajectories in the objective space obtained with the classical WS method and the proposed weight evolution scheme, together with the corresponding evolution of the weight parameter $q_1$ during the optimization. 
It can be observed that the adaptive weight evolution responds to the local structure of the objective landscape and consequently deviates from the trajectory produced by the standard WS method. This is also reflected in the evolution of the weight parameter $q_1$. In particular, after approximately $20$ iterations a clear deviation from the reference value becomes visible. This demonstrates that the proposed adaptive scheme is able to modify the effective weighting dynamically and thereby to explore regions of the Pareto frontier that would not be reached by a fixed weighted sum approach. 

\begin{figure}[ht]
    \centering
    \includegraphics[width=0.9\textwidth]{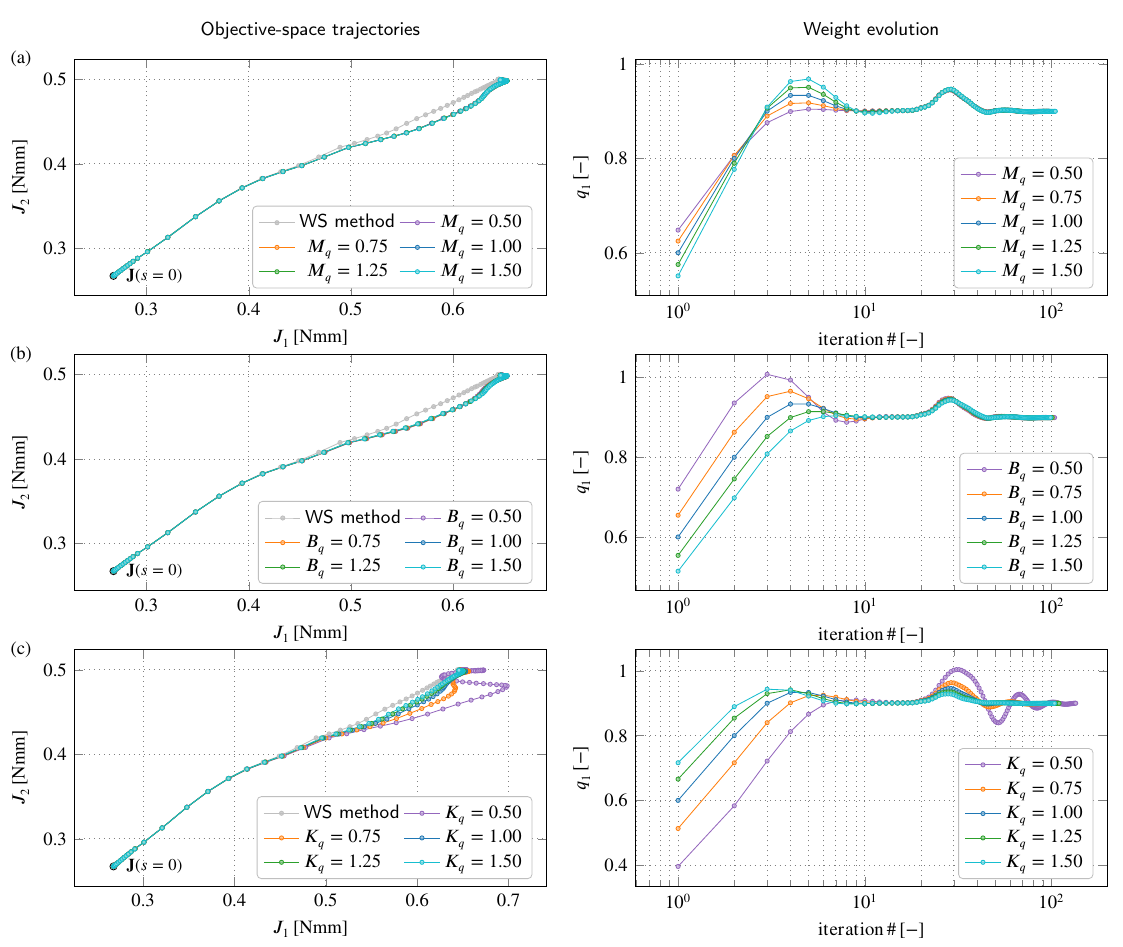}
    \caption{Effect of the numerical parameters on the objective-space trajectory and on the evolution of the weight parameter $q_1$: (a) effect of $M_q$, (b) effect of $B_q$, (c) effect of $K_q$.}
    \label{fig:expparameters}
\end{figure}

In the following, the effects of the numerical parameters on the resulting solution trajectory are investigated. As shown in Fig.~\ref{fig:expparameters}a, variations of $M_q$ have only a minor influence on the trajectory in the objective space. 
However, the effect becomes more apparent in the evolution of the weight parameter $q_1$, which is consistent with the previously derived mechanical analogy of the weight evolution scheme. As $M_q$ acts as an inertia term, larger values lead to a more pronounced overshooting behavior, whereas smaller values allow the weight parameter to adapt more rapidly to changes in the objective landscape. A similar behavior can be observed when varying the damping parameter $B_q$. Again, the trajectory itself is only weakly affected, while the main influence is visible in the behavior of $q_1$, as shown in Fig.~\ref{fig:expparameters}b. Increasing the damping suppresses oscillations and leads to a smoother convergence of the weight parameter $q_1$. A significantly stronger effect on the trajectory is obtained when varying $K_q$. With increasing $K_q$, the trajectory approaches that of the WS method, whereas smaller values of $K_q$ lead to a more sensitive response to variations in the objective landscape, thus controlling the scope of exploration as illustrated in Fig.~\ref{fig:expparameters}c. This finding is consistent with the previous analysis regarding the behavior of the weight evolution for $K_q\to\infty$, see Section~\ref{sec:methodology}.
\subsection{Influence of the initial conditions}
\begin{figure}[ht]
    \centering
    \includegraphics[width=0.9\textwidth]{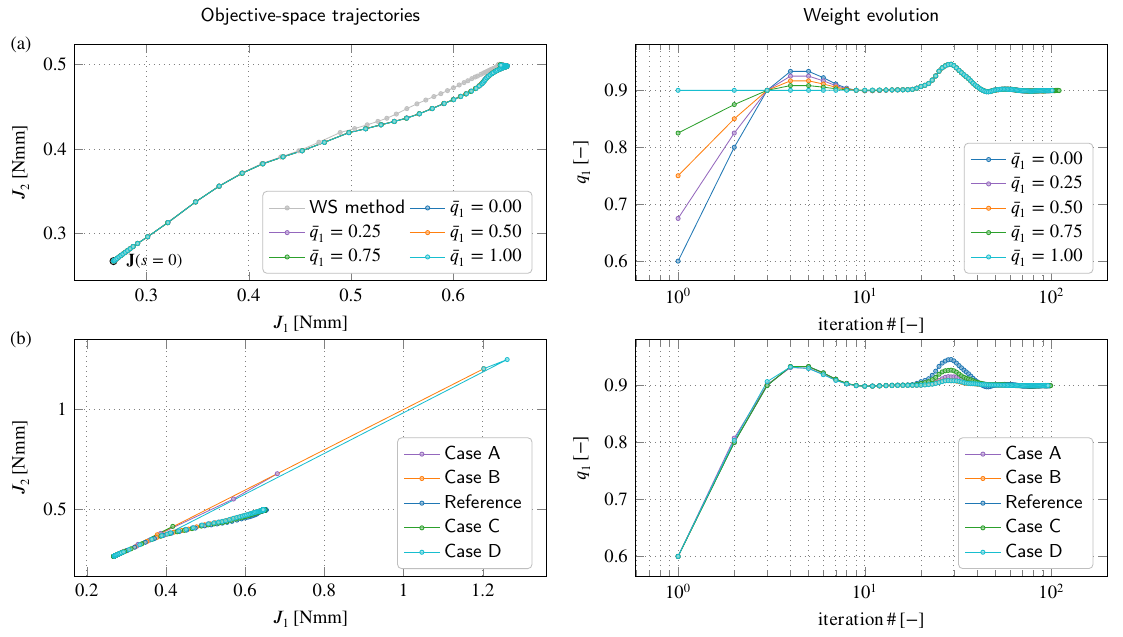}
    \caption{Effects of  the initial conditions of the coupled evolution problem: (a) effect of ratio $\bar{q}_1$, (b) effect of $\phi_{-1} $ and $\phi_0$.}
    \label{fig:expinitialconditions}
\end{figure}

\begin{figure}[ht]
    \centering
    \includegraphics[width=0.9\textwidth]{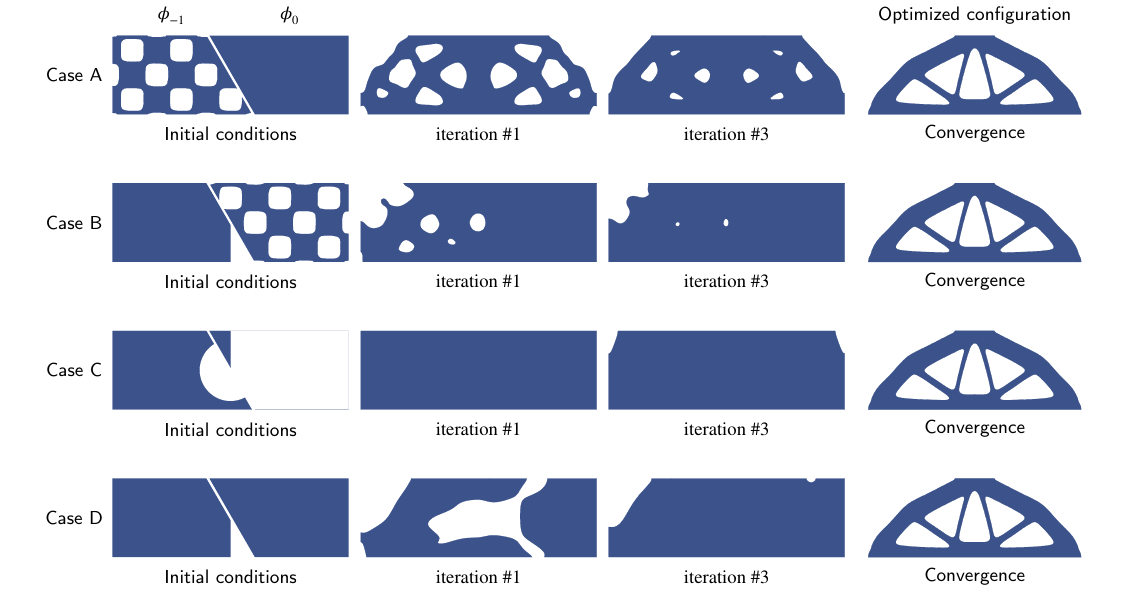}
    \caption{Convergence of topologies under different combinations of initial conditions $\phi_{-1}$ and $\phi_0$.}
    \label{fig:expinitialconditionstopo}
\end{figure}

Next, the influence of the initial conditions of the coupled evolution problem is investigated. First, the initial weight ratios $\bar{q}_1$ are varied. As illustrated in Fig.~\ref{fig:expinitialconditions}a, the resulting trajectories are nearly identical for all considered initial weight ratios. 
This indicates that the proposed evolution scheme is robust with respect to perturbations in the initial weights. In a second step, the sensitivity regarding the initial conditions of the level set function is analyzed by considering different combinations of $\phi_{-1}$ and $\phi_0$. The results show that perturbations in the topological layout are rapidly damped and that the optimization converges to almost identical topologies for all investigated configurations, see Fig.~\ref{fig:expinitialconditionstopo}. 
Furthermore, the corresponding objective-space trajectories shown in Fig.~\ref{fig:expinitialconditions}b reveal only minor deviations.

Small differences become noticeable after approximately $20$ iterations, where the algorithm interacts with the objective landscape in slightly different ways depending on the initial conditions. This behavior can be explained by the fact that the corresponding region of the landscape is approached from different directions. However, $q_1$ eventually converges to the same final state in all cases. Overall, the experiments indicate that the proposed coupled evolution scheme is robust with respect to variations in both the initial weights and the initial level set configurations.
\subsubsection{Adaptive Pareto frontier refinement using ASD}
\begin{figure}[ht]
    \centering
    \includegraphics[width=0.9\textwidth]{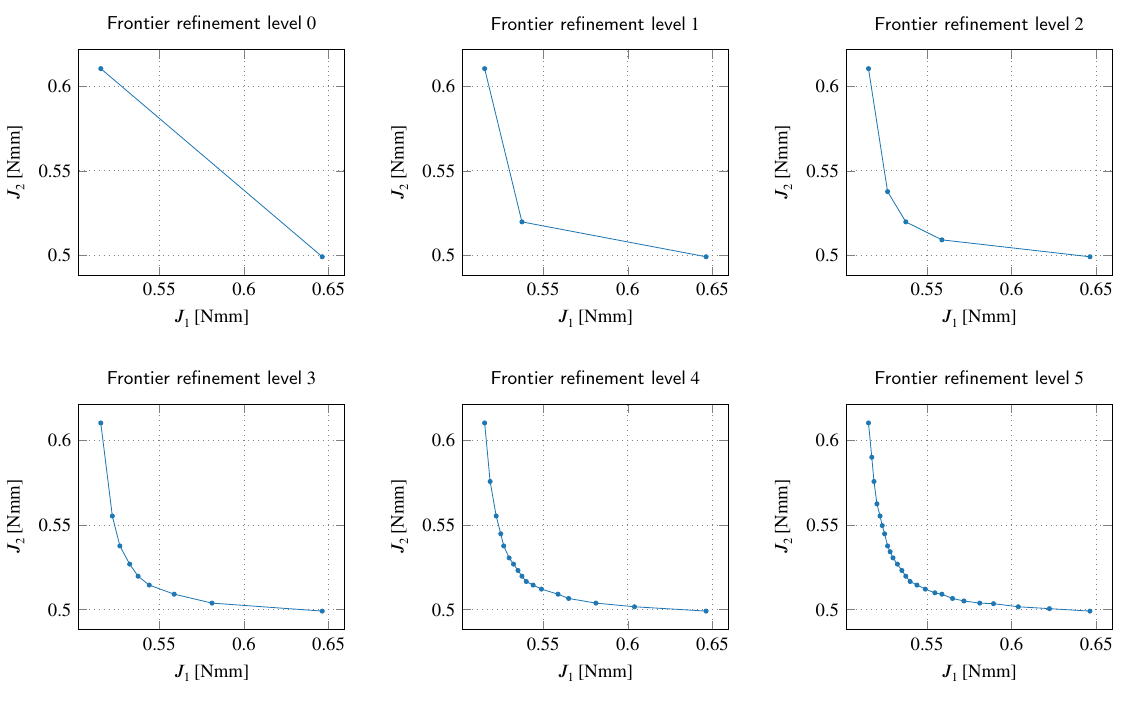}
    \caption{Adaptive Pareto frontier refinement of the bi-objective minimum mean compliance problem using the ASD method.}
    \label{fig:GRDfrontierrefinement}
\end{figure}

\begin{figure}[ht]
    \centering
    \includegraphics[width=0.9\textwidth]{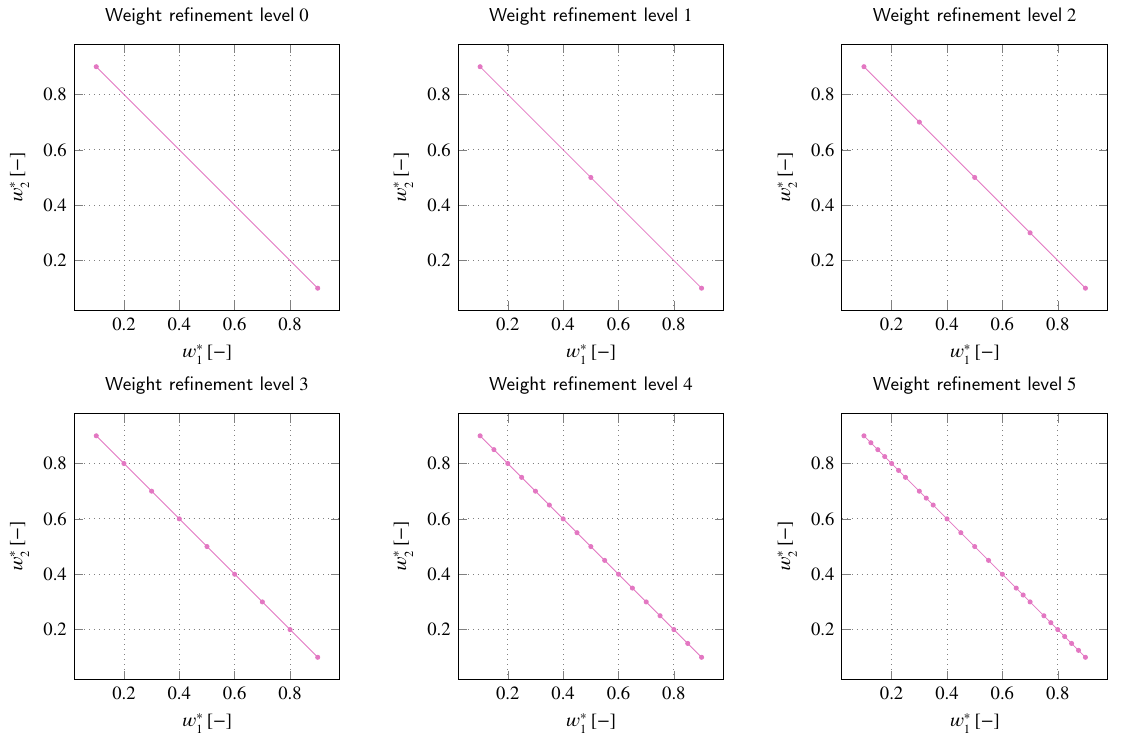}
    \caption{Weight refinement of the bi-objective minimum mean compliance problem.}
    \label{fig:GRDweightrefinement}
\end{figure}

\begin{figure}[ht]
    \centering
    \includegraphics[width=0.9\textwidth]{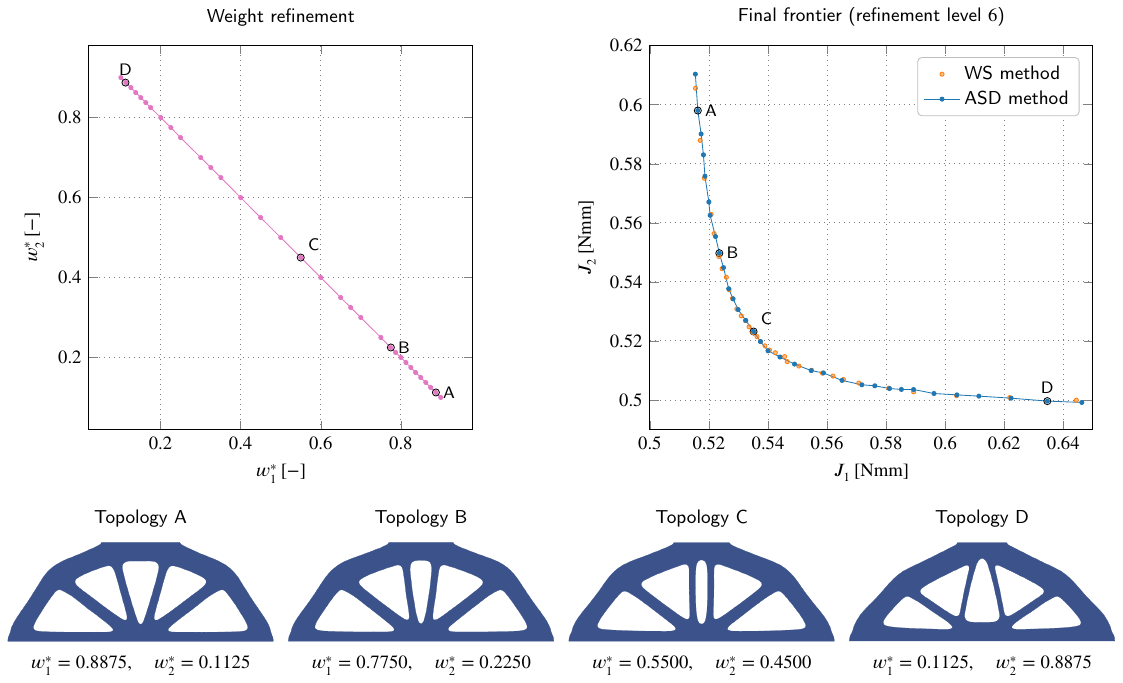}
    \caption{Final Pareto frontier approximation and selected topologies of the bi-objective minimum mean compliance problem.}
    \label{fig:GRDfrontierfinal}
\end{figure}

Finally, the performance of the proposed ASD method is demonstrated by adaptively refining the Pareto frontier using parameters $M_q=0.5$, $B_q=2$, and $K_q=1$. The initial weight discretization $\mathcal{W}_h^\ast=\{(0.9,0.1),(0.1,0.9)\}$ is chosen together with a threshold value of $l_\mathrm{s,max} = 0.04$. Figure~\ref{fig:GRDfrontierrefinement} shows the successive approximation of the Pareto frontier.
According to the principle of the method, the poor simplices are subdivided at each refinement level, which promotes a more uniform distribution of solutions along the frontier. The decomposition of the weight simplex is shown in Fig.~\ref{fig:GRDweightrefinement}.
The final approximation of the Pareto frontier is presented in Fig.~\ref{fig:GRDfrontierfinal}. Representative topologies along the frontier illustrate the structural trade-offs between the objectives, which become visible through noticeable rearrangements of the internal struts depending on the relative objective weighting.
For comparison, the results obtained using the WS method are also depicted. As expected, the WS method leads to the typical clustering of solutions in certain regions of the Pareto frontier.
\begin{figure}[ht]
    \centering
    \includegraphics[width=0.9\textwidth]{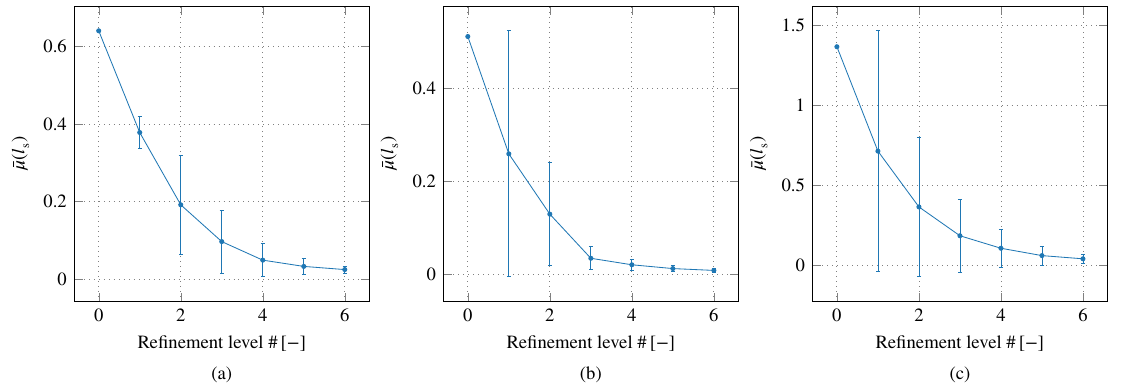}
    \caption{Convergence of the mean simplex length and its standard deviation: (a) bi-objective minimum mean compliance problem, (b) gripper design problem, (c) L-bracket problem.}
    \label{fig:simplexconvergence}
\end{figure}
The convergence behavior of the mean simplex edge length is shown in Fig.~\ref{fig:simplexconvergence}a. Both the mean value and the standard deviation decrease rapidly with increasing refinement level, indicating a progressive homogenization of the discretization. This confirms that the proposed ASD method provides a stable and systematic refinement of the Pareto frontier. All computed solution candidates were verified to be Pareto efficient in the final step of the procedure. Moreover, no points were removed when applying a point-distance tolerance of $1\times10^{-3}$. This results in a well-distributed approximation of the Pareto frontier without the clustering effects observed for the classical WS approach.
\subsection{Bi-objective gripper design}
\begin{figure}[ht]
    \centering
    \includegraphics[width=0.9\textwidth]{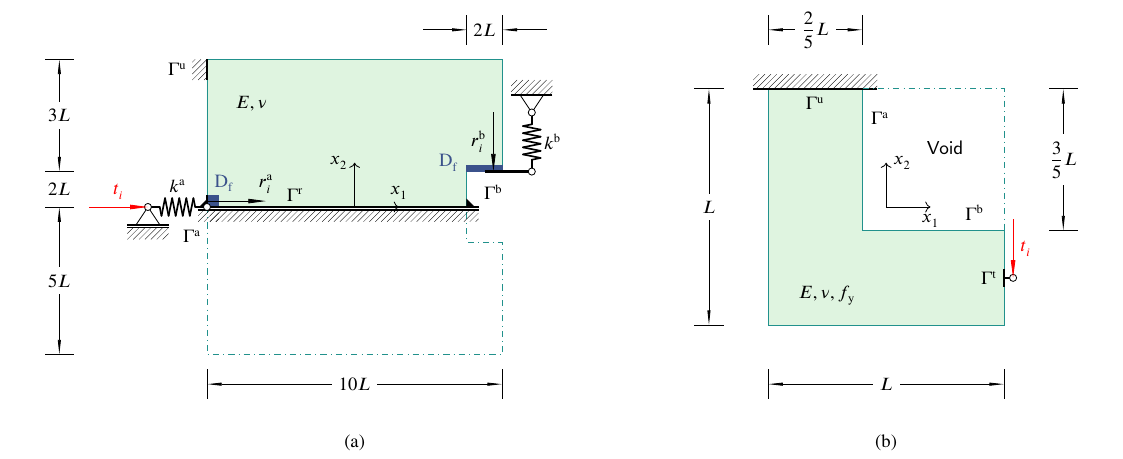}
    \caption{Bi-objective optimization problems: (a) gripper design problem, (b) L-bracket problem.}
    \label{fig:GRPLBRproblems}
\end{figure}

\begin{figure}[ht]
    \centering
    \includegraphics[width=0.9\textwidth]{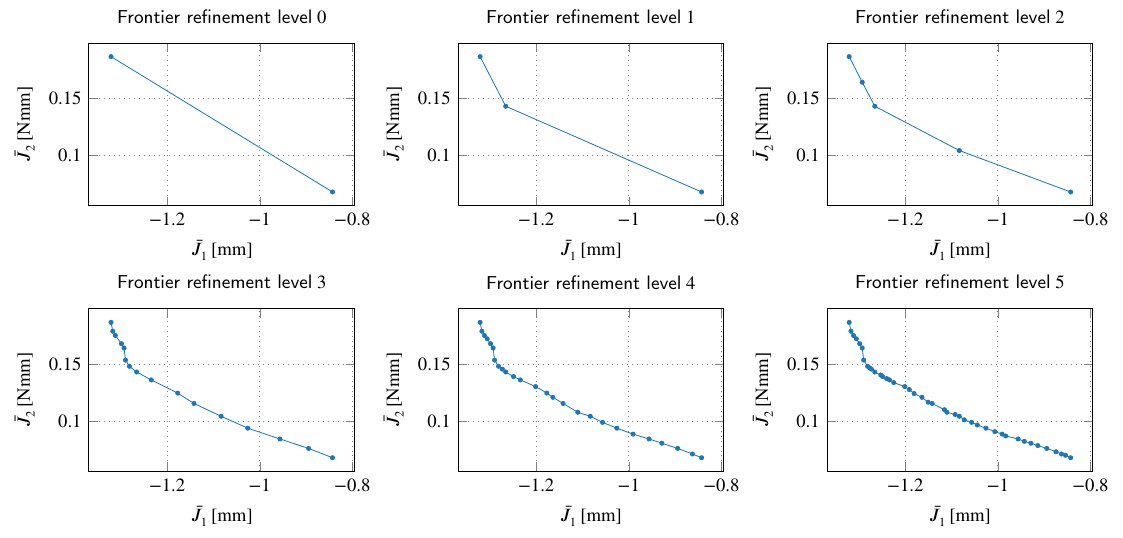}
    \caption{Adaptive Pareto frontier refinement for the gripper design problem; objective functionals scaled as $\bar{J}_1=J_1\times 10^{10}$ and $\bar{J}_2=J_2\times 10^{11}$.}
    \label{fig:GRPfrontierrefinement}
\end{figure}

\begin{figure}[ht]
    \centering
    \includegraphics[width=0.9\textwidth]{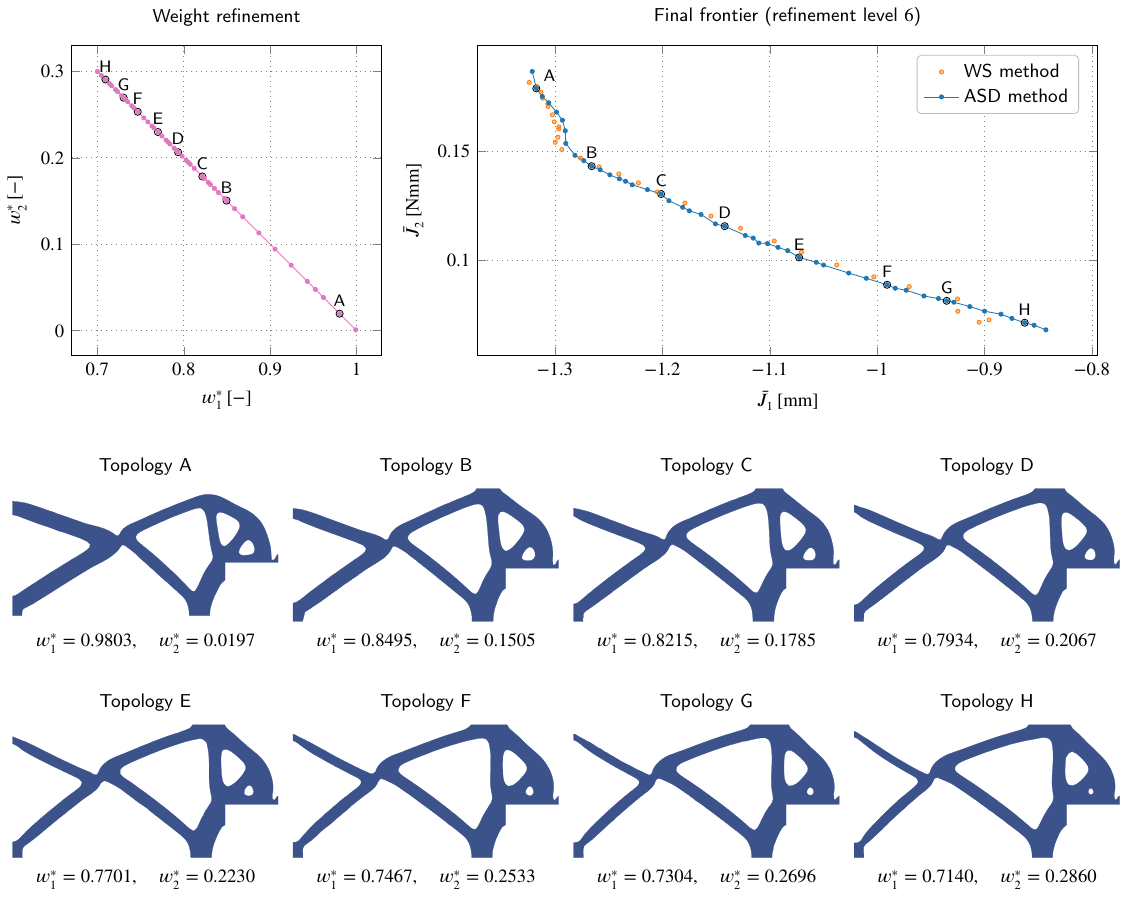}
    \caption{Final Pareto frontier approximation and selected topologies for the gripper design problem considering scaled objective functionals $\bar{J}_1=J_1\times 10^{10}$ and $\bar{J}_2=J_2\times 10^{11}$.}
    \label{fig:GRPfinalfrontier}
\end{figure}

Next, the proposed method is applied to a bi-objective compliant mechanism design problem ($\alpha\in\{1,2\}$), illustrated in Fig.~\ref{fig:GRPLBRproblems}a. 
The aim is to design a gripper mechanism that maximizes the output displacement while minimizing the strain energy. The material domain $\Omega$ is fixed along the boundary section $\Gamma^\mathrm{u}$ and embedded within the design domain $\mathrm{D}$, which additionally contains a fixed region of solid material $\mathrm{D}_\mathrm{f}$. Owing to the symmetry of the system, the boundary section $\Gamma^\mathrm{r}$ constrains the displacement in the normal direction. An input force $t_i$ is applied on boundary $\Gamma^\mathrm{a}$, while $\Gamma^\mathrm{b}$ corresponds to the output displacement region. The vectors $\mathbf{r}^\mathrm{a},\mathbf{r}^\mathrm{b}\in\mathbb{R}^{n\times 1}$ define the desired directions of motion at the respective boundaries. Artificial springs with stiffnesses $k^\mathrm{a}$ and $k^\mathrm{b}$ are introduced on both sections to generate reaction forces and enhance numerical stability \citep{SIGMUND1997}. Considering an additional volume constraint and neglecting body forces, the optimization problem is stated as
\begin{customopti}|s|
    {vec\quad inf}{\phi\in H^2(\mathscr{D})}{
    \mathbf{J}[\mathbf{u}_1,\mathbf{u}_2,\Theta] = \left(J_1[\mathbf{u}_1],J_2[\mathbf{u}_2,\Theta]\right)^\top}{}{}{}
    \addConstraint{G[\Theta] = \frac{1}{V_0}\int_{\mathrm{D}\setminus\mathrm{D}_\mathrm{f}}\Theta\,\mathrm{d}\Omega - V_\mathrm{f}}{\leq 0}{}{}
    \addConstraint{-\partial_j(\mathbb{C}_{ijkl}\varepsilon_{kl}(\mathbf{u}_\alpha))}{=0}{\quad\text{in}\quad\Omega}{}
    \addConstraint{u_{\alpha,i}}{=0}{\quad\text{on}\quad\Gamma^\mathrm{u}}{}
    \addConstraint{u_{\alpha,i} n_i}{=0}{\quad\text{on}\quad\Gamma^\mathrm{r}}{}
    \addConstraint{(\mathbb{C}_{ijkl}\varepsilon_{kl}(\mathbf{u}_\alpha))n_j}{=t_i- k^\mathrm{a}(r^\mathrm{a}_i r^\mathrm{a}_j){u}_{\alpha,j}}{\quad\text{on}\quad\Gamma^\mathrm{a}}{}
    \addConstraint{(\mathbb{C}_{ijkl}\varepsilon_{kl}(\mathbf{u}_\alpha))n_j}{=- k^\mathrm{b}(r^\mathrm{b}_i r^\mathrm{b}_j)u_{\alpha,j}}{\quad\text{on}\quad\Gamma^\mathrm{b}}{}
    \addConstraint{(\mathbb{C}_{ijkl}\varepsilon_{kl}(\mathbf{u}_\alpha))n_j}{=0}{\quad\text{on}\quad\partial\Omega\setminus\Gamma^\mathrm{g}}{}
    \label{eq:gripperproblem}
\end{customopti}  
with objective functionals
\begin{align*}
    J_1[\mathbf{u}_1] &= -\int_{\Gamma^\mathrm{b}} r^\mathrm{b}_i u_{1,i} \,\mathrm{d}\Gamma, \\ 
    J_2[\mathbf{u}_2,\Theta] &= \frac{1}{2}\int_{\mathrm{D}\setminus\mathrm{D}_\mathrm{f}} \tau(\Theta)\mathbb{C}_{ijkl}\varepsilon_{kl}(\mathbf{u}_2)\varepsilon_{ij}(\mathbf{u}_2) \,\mathrm{d}\Omega + \frac{1}{2}\int_{\mathrm{D}_\mathrm{f}}\mathbb{C}_{ijkl}\varepsilon_{kl}(\mathbf{u}_2)\varepsilon_{ij}(\mathbf{u}_2)\,\mathrm{d}\Omega
\end{align*}
and summarized boundary sections $\Gamma^\mathrm{g} \vcentcolon= \Gamma^\mathrm{u}\cup\Gamma^\mathrm{r}\cup\Gamma^\mathrm{a}\cup\Gamma^\mathrm{b}$. From the boundary conditions in Problem~(\ref{eq:gripperproblem}), the governing equations in weak form are obtained
\begin{align}
    R_\alpha:\quad\int_{\mathrm{D}\setminus\mathrm{D}_\mathrm{f}}\tau(\Theta)\mathbb{C}_{ijkl}\varepsilon_{kl}(\mathbf{u}_\alpha)\varepsilon_{ij}(\mathbf{v}_\alpha)\,\mathrm{d}\Omega + \int_{\mathrm{D}_\mathrm{f}}\mathbb{C}_{ijkl}\varepsilon_{kl}(\mathbf{u}_\alpha)\varepsilon_{ij}(\mathbf{v}_\alpha)\,\mathrm{d}\Omega + \dots \notag\\
    \dots + \int_{\Gamma^\mathrm{b}}k^\mathrm{b}(r^{\mathrm{b}}_i r^{\mathrm{b}}_j)u_{\alpha,j} v_{\alpha,i}\,\mathrm{d}\Gamma = \int_{\Gamma^\mathrm{a}}(t_i - k^\mathrm{a}(r^{\mathrm{a}}_i r^{\mathrm{a}}_j)u_{\alpha,j})v_{\alpha,i}\,\mathrm{d}\Gamma
    \label{eq:GRPgoverning}
\end{align}
and the adjoint equations follow as
\begin{align}
    \mathcal{A}_\alpha:\quad \left\langle \frac{w_\alpha}{J_\alpha^\ast}\frac{\partial J_\alpha}{\partial\mathbf{u}_\alpha},\delta\mathbf{u}_\alpha\right\rangle &= \int_{\mathrm{D}\setminus\mathrm{D}_\mathrm{f}}\tau(\Theta)\mathbb{C}_{ijkl}\varepsilon_{kl}(\delta\mathbf{u}_\alpha)\varepsilon_{ij}(\mathbf{v}_\alpha)\,\mathrm{d}\Omega + \dots \notag\\
    &\dots + \int_{\mathrm{D}_\mathrm{f}}\mathbb{C}_{ijkl}\varepsilon_{kl}(\delta\mathbf{u}_\alpha)\varepsilon_{ij}(\mathbf{v}_\alpha)\,\mathrm{d}\Omega + \dots \notag\\
    &\dots + \int_{\Gamma^\mathrm{b}}k^\mathrm{b}(r^\mathrm{b}_i r^\mathrm{b}_j)\delta u_{\alpha,i}v_{\alpha,i}\,\mathrm{d}\Gamma + \dots \notag\\
    &\dots + \int_{\Gamma^\mathrm{a}}k^\mathrm{a}(r^{\mathrm{a}}_i r^{\mathrm{a}}_j)\delta u_{\alpha,j}v_{\alpha,i}\,\mathrm{d}\Gamma
    \label{eq:GRPadjoint}
\end{align}
where
\begin{align*}
    \left\langle \frac{w_1}{J_1^\ast}\frac{\partial J_1}{\partial\mathbf{u}_1},\delta\mathbf{u}_1\right\rangle &= -\int_{\Gamma^\mathrm{b}}\frac{w_1}{J_1^\ast}r^{\mathrm{b}}_i\delta u_{1,i}\,\mathrm{d}\Gamma, \\
    \left\langle \frac{w_2}{J_2^\ast}\frac{\partial J_2}{\partial\mathbf{u}_2},\delta\mathbf{u}_2\right\rangle &= \int_{\mathrm{D}\setminus\mathrm{D}_\mathrm{f}}\frac{w_2}{J_2^\ast}\tau(\Theta)\mathbb{C}_{ijkl}\varepsilon_{kl}(\delta\mathbf{u}_2)\varepsilon_{ij}(\mathbf{u}_2)\,\mathrm{d}\Omega + \dots\\ 
    &\dots + \int_{\mathrm{D}_\mathrm{f}}\frac{w_2}{J_2^\ast}\mathbb{C}_{ijkl}\varepsilon_{kl}(\delta\mathbf{u}_2)\varepsilon_{ij}(\mathbf{u}_2)\,\mathrm{d}\Omega. 
\end{align*}
The partial derivatives of the objective functionals with respect to $\Theta$ are given by
\begin{align*}
    \frac{w_1}{J_1^\ast}\frac{\partial J_1}{\partial\Theta} &= 0, \\
    \frac{w_2}{J_2^\ast}\frac{\partial J_2}{\partial\Theta} &= \frac{w_2}{2 J_2^\ast}\frac{\partial\tau}{\partial\Theta}\mathbb{C}_{ijkl}\varepsilon_{kl}(\mathbf{u}_2)\varepsilon_{ij}(\mathbf{u}_2).
\end{align*}
The derivative of the constraint functional corresponds to Eq.~(\ref{eq:GRDpartialderivativeconstraint}), while the derivatives of the governing equations are identical to those given in Eq.~(\ref{eq:GRDgoverningequations}). Subsequently, one obtains the individual perturbation terms
\begin{align}
    f_1 &= \frac{\lambda}{2V_0} - \frac{\partial\tau}{\partial\Theta}\mathbb{C}_{ijkl}\varepsilon_{kl}(\mathbf{u}_1)\varepsilon_{ij}(\mathbf{v}_1),\label{eq:GRPperturbationterm1}\\
    f_2 &= \frac{\lambda}{2V_0} + \frac{w_2}{2J_2^\ast}\frac{\partial\tau}{\partial\Theta}\mathbb{C}_{ijkl}\varepsilon_{kl}(\mathbf{u}_2)\varepsilon_{ij}(\mathbf{u}_2) - \frac{\partial\tau}{\partial\Theta}\mathbb{C}_{ijkl}\varepsilon_{kl}(\mathbf{u}_2)\varepsilon_{ij}(\mathbf{v}_2)\label{eq:GRPperturbationterm2}.
\end{align}
Summation of Eq.~(\ref{eq:GRPperturbationterm1}) and Eq.~(\ref{eq:GRPperturbationterm2}) finally yields the overall perturbation term:
\begin{align}
    F = \frac{\lambda}{V_0} + \frac{\partial\tau}{\partial\Theta}\mathbb{C}_{ijkl}\left[- \varepsilon_{kl}(\mathbf{u}_1)\varepsilon_{ij}(\mathbf{v}_1) + \frac{w_2}{2J_2^\ast}\varepsilon_{kl}(\mathbf{u}_2)\varepsilon_{ij}(\mathbf{u}_2) - \varepsilon_{kl}(\mathbf{u}_2)\varepsilon_{ij}(\mathbf{v}_2)\right]
    \label{eq:GRPperturbationterm}.
\end{align}
The normalization factors correspond to those of Eq.~(\ref{eq:GRDnormalizationfactor}), thus the evolution of the level set function is fully determined by Eqs.~(\ref{eq:GRPgoverning}), (\ref{eq:GRPadjoint}) and (\ref{eq:GRPperturbationterm}). For the numerical experiments, the artificial spring stiffnesses are set to $k^\mathrm{a} = 1\times 10^5\,\mathrm{N}/\mathrm{mm}$ and $k^\mathrm{b} = 1\times 10^3\,\mathrm{N}/\mathrm{mm}$ and the direction vectors to $\mathbf{r}^\mathrm{a}=[1,0]^\top$ and $\mathbf{r}^\mathrm{b}=[0,-1]^\top$. The volume fraction is set to $V_\mathrm{f} = 0.30$. The wave propagation components are specified as $C_{11} = C_{22} = 0.011$, and a damping coefficient of $B = 0.1$ is considered. The parameters governing the weight evolution are chosen as $M_q = 8$, $B_q = 14$, and $K_q = 10$, while the maximum simplex length is set to $l_\mathrm{s,max} = 0.01$. The initial reference weight discretization is chosen as $\mathcal{W}_h^\ast=\{(0.999,0.001),(0.70,0.30)\}$. In this example, we again observe that the ASD method is able to gradually refine the Pareto frontier, as illustrated in Fig.~\ref{fig:GRPfrontierrefinement}. 
In the early iterations, the simplex lengths exhibit significant variability. However, this difference is quickly reduced through adaptive refinement, see Fig~\ref{fig:simplexconvergence}b. Figure~\ref{fig:GRPfrontierrefinement} demonstrates how the frontier is progressively stabilized and finely resolved over the iterations. The final Pareto frontier, along with selected example topologies, is shown in Fig.~\ref{fig:GRPfinalfrontier}.
For comparison, the results obtained with the conventional WS method are also included. This highlights that the ASD method produces a more stable approximation of the Pareto frontier, whereas the WS method exhibits stronger divergence, for instance in the region between topologies G and H. Moreover, a local non-convexity is observed between topologies A and B, which is consistently captured and refined by the ASD method. A structural analysis of the topologies indicates that, for maximum output displacement, a more robust geometry of the bracket on the input side is preferred. As the priority placed on minimizing strain energy increases, the design of the more robust geometry shifts toward the output side. In the final step, nearly coincident points were removed using a distance tolerance of $6\times 10^{-3}$. All remaining solution candidates were verified to be Pareto efficient.
\subsection{Bi-objective volume and strain energy minimization}
\begin{figure}[ht]
    \centering
    \includegraphics[width=0.9\textwidth]{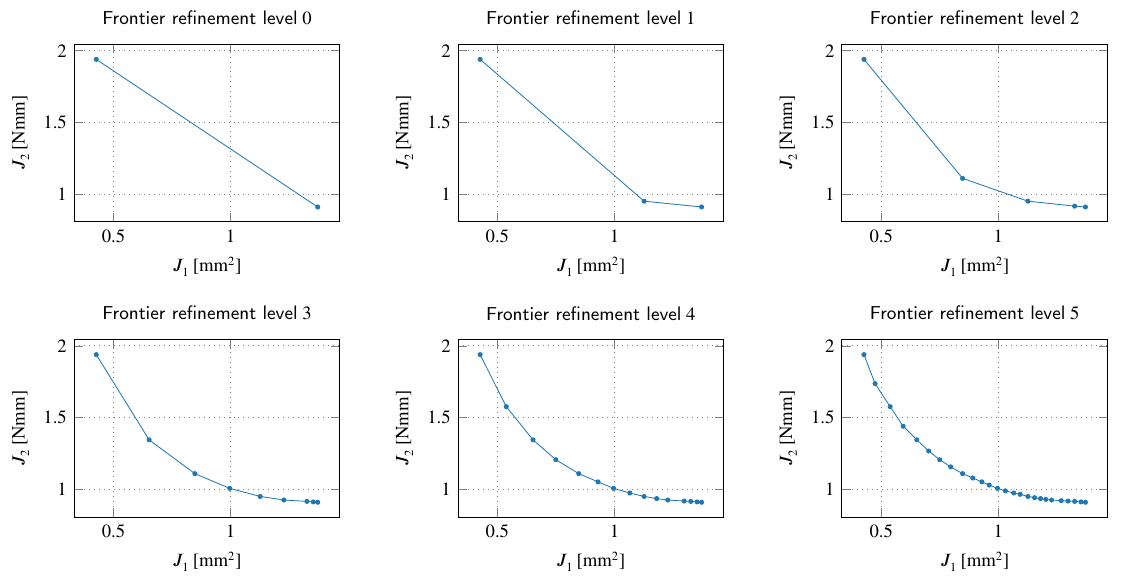}
    \caption{Adaptive Pareto frontier refinement for the L-bracket problem.}
    \label{fig:LBRfrontierrefinement}
\end{figure}

\begin{figure}[ht]
    \centering
    \includegraphics[width=0.9\textwidth]{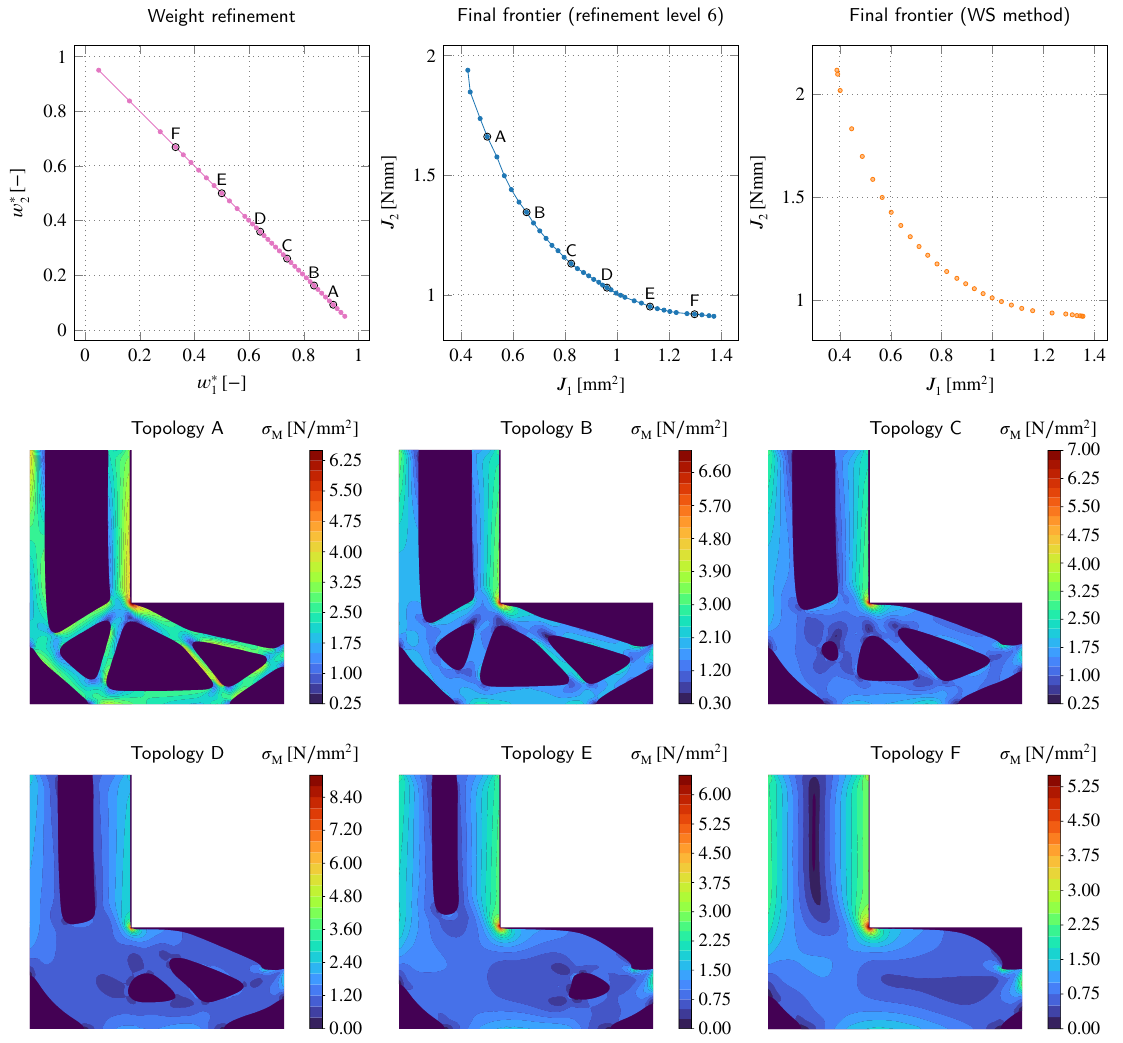}
    \caption{Final Pareto frontier approximation and selected topologies for the L-bracket problem.}
    \label{fig:LBRfrontierfinal}
\end{figure}

We now investigate a fundamental bi-objective trade-off in lightweight structural optimization, balancing material efficiency against structural stiffness ($\alpha\in\{1,2\}$) applied to the L-bracket problem depicted in Fig.~\ref{fig:GRPLBRproblems}b. To account for stress limitations, an additional global constraint based on a $p$-norm aggregation of the von Mises stress is introduced \citep{LE2010,PICELLI2018}. Consider a material domain $\Omega$ fixed at $\Gamma^\mathrm{u}$ and loaded by a traction vector $t_i$ at $\Gamma^\mathrm{t}$. To represent adjacent void regions, Dirichlet boundary conditions $\phi = -1$ are imposed along $\Gamma^\mathrm{a}$ and $\Gamma^\mathrm{b}$. The resulting bi-objective optimization problem is then formulated as
\begin{customopti}|s|
    {vec\quad inf}{\phi\in H^2(\mathscr{D})}{
    \mathbf{J}[\mathbf{u}_1,\mathbf{u}_2, \Theta] = \left(J_1[\Theta],J_2[\mathbf{u}_2, \Theta]\right)^\top}{}{}{}
    \addConstraint{G_\alpha[\mathbf{u}_\alpha,\Theta] = \frac{1}{V_0}S_\sigma^{1/p}[\mathbf{u}_\alpha, \Theta;p] - \bar{\sigma}}{\leq 0}{}{}
    \addConstraint{-\partial_j(\mathbb{C}_{ijkl}\varepsilon_{kl}(\mathbf{u}_\alpha))}{=0}{\quad\text{in}\quad\Omega}{}
    \addConstraint{u_{\alpha,i}}{=0}{\quad\text{on}\quad\Gamma^\mathrm{u}}{}
    \addConstraint{(\mathbb{C}_{ijkl}\varepsilon_{kl}(\mathbf{u}_\alpha))n_j}{=t_i}{\quad\text{on}\quad\Gamma^\mathrm{t}}{}
    \addConstraint{(\mathbb{C}_{ijkl}\varepsilon_{kl}(\mathbf{u}_\alpha))n_j}{=0}{\quad\text{on}\quad\partial\Omega\setminus(\Gamma^\mathrm{u}\cup\Gamma^\mathrm{t})}{}
    \label{eq:selfweightproblem}
\end{customopti}  
with objective functionals
\begin{align*}
    J_1[\Theta] &= \int_\mathrm{D}\Theta\,\mathrm{d}\Omega,\\
    J_2[\mathbf{u}_2,\Theta] &= \frac{1}{2}\int_\mathrm{D}\tau(\Theta)\mathbb{C}_{ijkl}\varepsilon_{kl}(\mathbf{u}_2)\varepsilon_{ij}(\mathbf{u}_2)\,\mathrm{d}\Omega
\end{align*}
and aggregated stress functional
\begin{align*}
    S_\sigma^{1/p}[\mathbf{u}_\alpha,\Theta;p] = \left(\int_\mathrm{D}\left(\frac{\sigma_\mathrm{M}(\mathbf{u}_\alpha)}{f_\mathrm{y}}\right)^p\tau(\Theta)\,\mathrm{d}\Omega\right)^{1/p}.
\end{align*}
Here, $f_\mathrm{y}$ represents the maximum admissible stress while $p$ controls the penalty of stress concentrations and $\bar{\sigma}$ denotes the limit value for the global stress aggregation. Using the stress deviator tensor
\begin{align*}
    s_{ij}(\mathbf{u}_\alpha) = \sigma_{ij}(\mathbf{u}_\alpha) - \frac{1}{3}\sigma_{kk}(\mathbf{u}_\alpha)\delta_{ij}
\end{align*}
the von Mises stress is expressed by the relation
\begin{align*}
    \sigma_\mathrm{M}(\mathbf{u}_\alpha) = \sqrt{\frac{3}{2}s_{ij}(\mathbf{u}_\alpha)s_{ij}(\mathbf{u}_\alpha)}.
\end{align*}
First, we derive the governing equations
\begin{align}
    R_\alpha:\quad\int_\mathrm{D}\tau(\Theta)\mathbb{C}_{ijkl}\varepsilon_{kl}(\mathbf{u}_\alpha)\varepsilon_{ij}(\mathbf{v}_\alpha)\,\mathrm{d}\Omega - \int_{\Gamma^\mathrm{t}}t_i v_{\alpha,i}\,\mathrm{d}\Gamma = 0,\quad\mathbf{v}_\alpha\in\mathcal{V}_\alpha
    \label{eq:LBRgoverning}
\end{align}
from the boundary conditions in Problem~(\ref{eq:selfweightproblem}). The adjoint equations are obtained as
\begin{align}
    \mathcal{A}_\alpha:\quad \left\langle\lambda_\alpha\frac{\partial G}{\partial\mathbf{u}_\alpha},\delta\mathbf{u}_\alpha\right\rangle + \left\langle\frac{w_\alpha}{J_\alpha^\ast}\frac{\partial J_\alpha}{\partial\mathbf{u}_\alpha},\delta\mathbf{u}_\alpha\right\rangle = \dots\notag\\
    \dots \int_\mathrm{D}\tau(\Theta)\mathbb{C}_{ijkl}\varepsilon_{kl}(\delta\mathbf{u}_\alpha)\varepsilon_{ij}(\mathbf{v}_\alpha)\,\mathrm{d}\Omega,\quad\delta\mathbf{u}_\alpha\in\mathcal{V}_\alpha
    \label{eq:LBRadjoint}
\end{align}
with the individual functional derivatives
\begin{align*}
    \left\langle\frac{w_1}{J_1^\ast}\frac{\partial J_1}{\partial\mathbf{u}_1},\delta\mathbf{u}_1\right\rangle &= 0, \\
    \left\langle\frac{w_2}{J_2^\ast}\frac{\partial J_2}{\partial\mathbf{u}_2},\delta\mathbf{u}_2\right\rangle &= \int_{\mathrm{D}} \frac{w_2}{J_2^\ast}\tau(\Theta)\mathbb{C}_{ijkl}\varepsilon_{kl}(\delta\mathbf{u}_2)\varepsilon_{ij}(\mathbf{u}_2) \,\mathrm{d}\Omega
\end{align*}
and
\begin{align*}
    \left\langle\lambda_\alpha\frac{\partial G_\alpha}{\partial\mathbf{u}_\alpha},\delta\mathbf{u}_\alpha\right\rangle = \frac{\lambda_\alpha}{V_0} S_\sigma^{1/p-1}\int_{\mathrm{D}}\left(\frac{\sigma_\mathrm{M}(\mathbf{u}_\alpha)}{f_\mathrm{y}}\right)^{p-1}\frac{3\tau(\Theta)}{2f_\mathrm{y}\sigma_\mathrm{M}(\mathbf{u}_\alpha)}s_{ij}(\mathbf{u}_\alpha)s_{ij}(\delta\mathbf{u}_\alpha)\,\mathrm{d}\Omega.
\end{align*}
The derivatives of the objective functionals with respect to $\Theta$ read
\begin{align*}
    \frac{w_1}{J_1^\ast}\frac{\partial J_1}{\partial\Theta} &= \frac{w_1}{J_1^\ast}, \\
    \frac{w_2}{J_2^\ast}\frac{\partial J_2}{\partial\Theta} &= \frac{w_2}{2J_2^\ast}\frac{\partial\tau}{\partial\Theta}\mathbb{C}_{ijkl}\varepsilon_{kl}(\mathbf{u}_2)\varepsilon_{ij}(\mathbf{u}_2).
\end{align*}
The derivative of the constraint functional for each problem with respect to $\Theta$ is given by
\begin{align*}
    \frac{\partial G_\alpha}{\partial\Theta}
    =
    \frac{\lambda_\alpha}{p V_0}
    S_{\sigma}^{1/p-1}[\mathbf{u}_\alpha,\Theta;p]
    \left(\frac{\sigma_\mathrm{M}(\mathbf{u}_\alpha)}{f_\mathrm{y}}\right)^p
    \frac{\partial\tau}{\partial\Theta}
\end{align*}
and the derivatives of the governing equations again correspond to those given in Eq.~(\ref{eq:GRDpartialderivativegoverning}). Thus, the overall perturbation term follows as
\begin{align}
    F = f_1 + f_2
    \label{eq:LBRperturbation}
\end{align}
with the respective contributions
\begin{align*}
    f_1 &= \frac{\lambda_1}{pV_0}S_\sigma^{1/p-1}[\mathbf{u}_1,\Theta;p]\left(\frac{\sigma_\mathrm{M}(\mathbf{u}_1)}{f_\mathrm{y}}\right)^p\frac{\partial\tau}{\partial\Theta} + \frac{w_1}{J_1^\ast} - \frac{\partial\tau}{\partial\Theta}\mathbb{C}_{ijkl}\varepsilon_{kl}(\mathbf{u}_1)\varepsilon_{ij}(\mathbf{v}_1),\\
    f_2 &= \frac{\lambda_2}{pV_0}S_\sigma^{1/p-1}[\mathbf{u}_2,\Theta;p]\left(\frac{\sigma_\mathrm{M}(\mathbf{u}_2)}{f_\mathrm{y}}\right)^p\frac{\partial\tau}{\partial\Theta} + \frac{w_2}{2J_2^\ast}\frac{\partial\tau}{\partial\Theta}\mathbb{C}_{ijkl}\varepsilon_{kl}(\mathbf{u}_2)\varepsilon_{ij}(\mathbf{u}_2) - \dots\\
    &\dots - \frac{\partial\tau}{\partial\Theta}\mathbb{C}_{ijkl}\varepsilon_{kl}(\mathbf{u}_2)\varepsilon_{ij}(\mathbf{v}_2).
\end{align*}
In the context of stress-constrained optimization, it is well known that structural modifications, such as phase separation or the nucleation of new holes, can lead to stress concentrations and thus high-frequency components in the perturbation term causing numerical instabilities. To mitigate these effects, we apply a nonlinear scaling following \citep{OELLERICH2025,EMMENDOERFER2016} and regularize the perturbation term $F$ from Eq.~(\ref{eq:LBRperturbation}) using a Helmholtz-type filtering:
\begin{align}
    \int_\mathrm{D}(\eta\partial_i \bar{F}\partial_i \delta\bar{F} + \bar{F}\delta\bar{F})\,\mathrm{d}\Omega 
    = \int_\mathrm{D} \frac{1}{\gamma}\mathrm{arsinh}(\gamma F)\delta \bar{F}\,\mathrm{d}\Omega,
    \label{eq:helmholtz}
\end{align}
where $\delta\bar{F}\in H^1(\mathrm{D})$ is a suitable test function and $\bar{F}$ denotes the regularized perturbation. Here, parameter $\eta > 0$ is the filtering parameter and $\gamma > 0$ controls the nonlinearity of the scaling. In this example, the normalization factors are determined as
\begin{align*}
    C_\alpha^s = \frac{1}{w_\alpha V_0}\int_\mathrm{D}\left\Vert\frac{\partial R_\alpha}{\partial\Theta} - \frac{w_\alpha}{J_\alpha^\ast}\frac{\partial J_\alpha}{\partial\Theta}\right\Vert_{L^2(\mathrm{D})}\,\mathrm{d}\Omega.
\end{align*}
The evolution problem of the level set function is now completely defined by Eqs.~(\ref{eq:LBRgoverning}), (\ref{eq:LBRadjoint}) and Eqs.~(\ref{eq:LBRperturbation}) and (\ref{eq:helmholtz}). For the numerical experiments, the fictitious material parameters are set to $C_{11}=C_{22}=0.015$ and the damping coefficient to $B=0.6$. For the stress aggregation, a penalization factor of $p=5$ and a yield stress of $f_\mathrm{y}=42\,\mathrm{N}/\mathrm{mm}^2$ are assumed, while the maximum admissible aggregated stress is limited to $\bar{\sigma}=0.05$. The parameters governing the weight evolution are chosen as $M_q=6$, $B_q=12$, and $K_q=10$, and the maximum simplex length is set to $l_\mathrm{s,max}=0.04$. For the regularization of the perturbation term, the values $\eta=1\times10^{-4}$ and $\gamma=2$ are used. The initial reference weight discretization is defined as $\mathcal{W}_h^\ast=\{(0.05,0.95),(0.95,0.05)\}$. This example further demonstrates the gradual refinement of the Pareto frontier achieved by the ASD method, see Fig.~\ref{fig:LBRfrontierrefinement}.
The resulting Pareto frontier, together with selected example topologies, is presented in Fig.~\ref{fig:LBRfrontierfinal}. 
As expected, increasing structural stiffness results in designs with higher material usage, whereas lighter configurations exhibit more slender load paths. The conventional WS method again shows clustering of solutions along the frontier. The corresponding convergence behavior of the ASD method is illustrated in Fig.~\ref{fig:simplexconvergence}c, where both the average simplex length and its standard deviation decrease rapidly over the iterations. In the final step, no points were removed for a point-distance tolerance of $1\times 10^{-3}$, and all solutions were verified to be Pareto efficient.
\subsection{Tri-objective minimum mean compliance}
\begin{figure}[ht]
    \centering
    \includegraphics[width=0.9\textwidth]{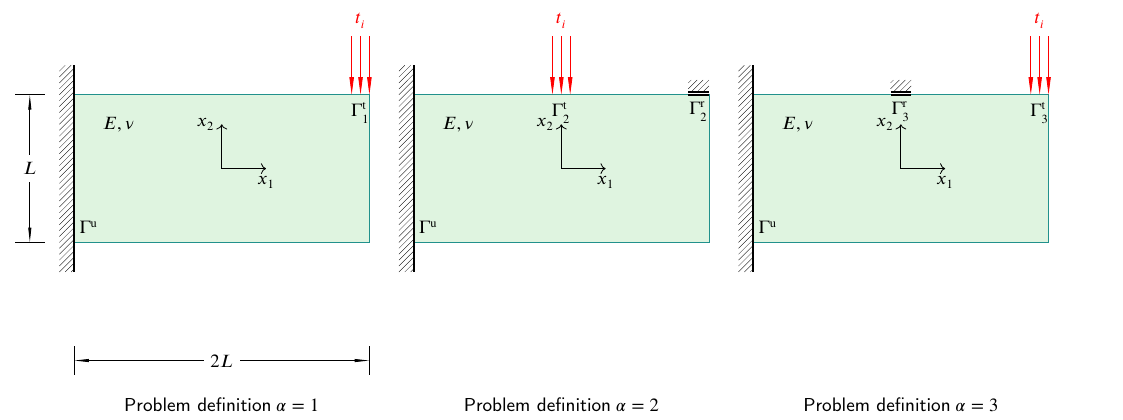}
    \caption{Tri-objective minimum mean compliance problem: clamped girder under different loading and supporting conditions.}
    \label{fig:problemCLB}
\end{figure}

\begin{figure}[ht]
    \centering
    \includegraphics[width=0.9\textwidth]{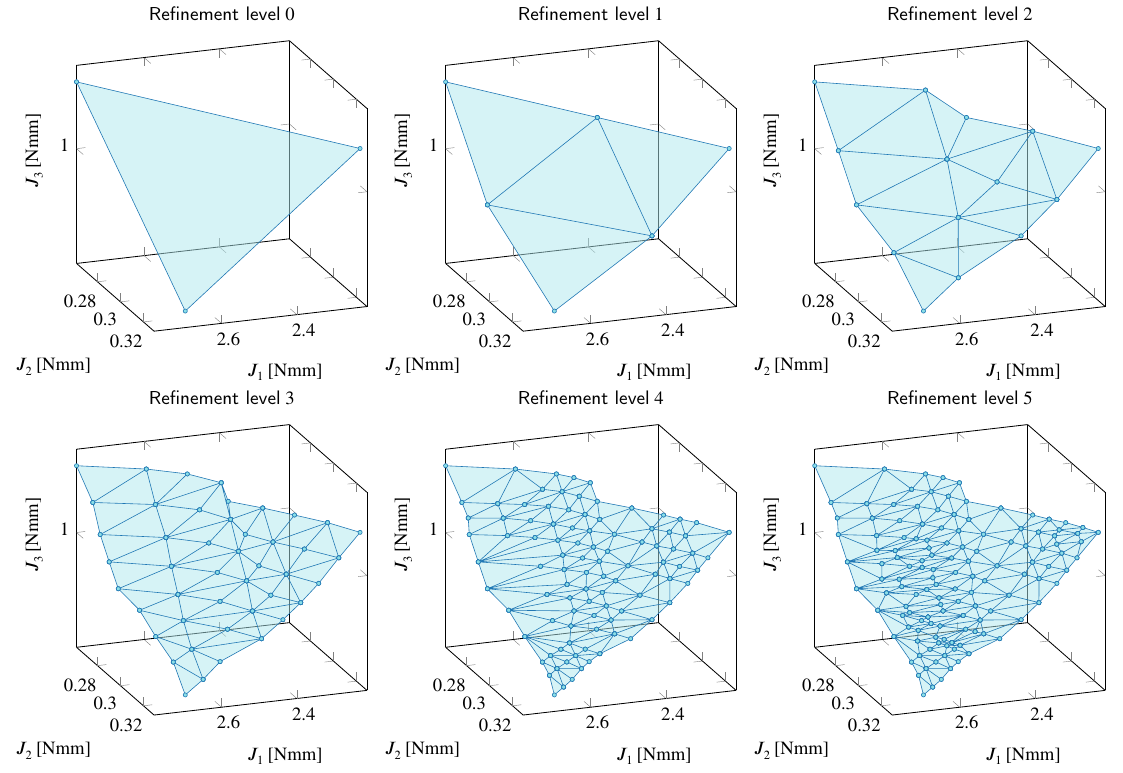}
    \caption{Adaptive Pareto frontier refinement of the tri-objective minimum mean compliance problem using the ASD method.}
    \label{fig:CLBfrontierrefinement}
\end{figure}

\begin{figure}[ht]
    \centering
    \includegraphics[width=0.9\textwidth]{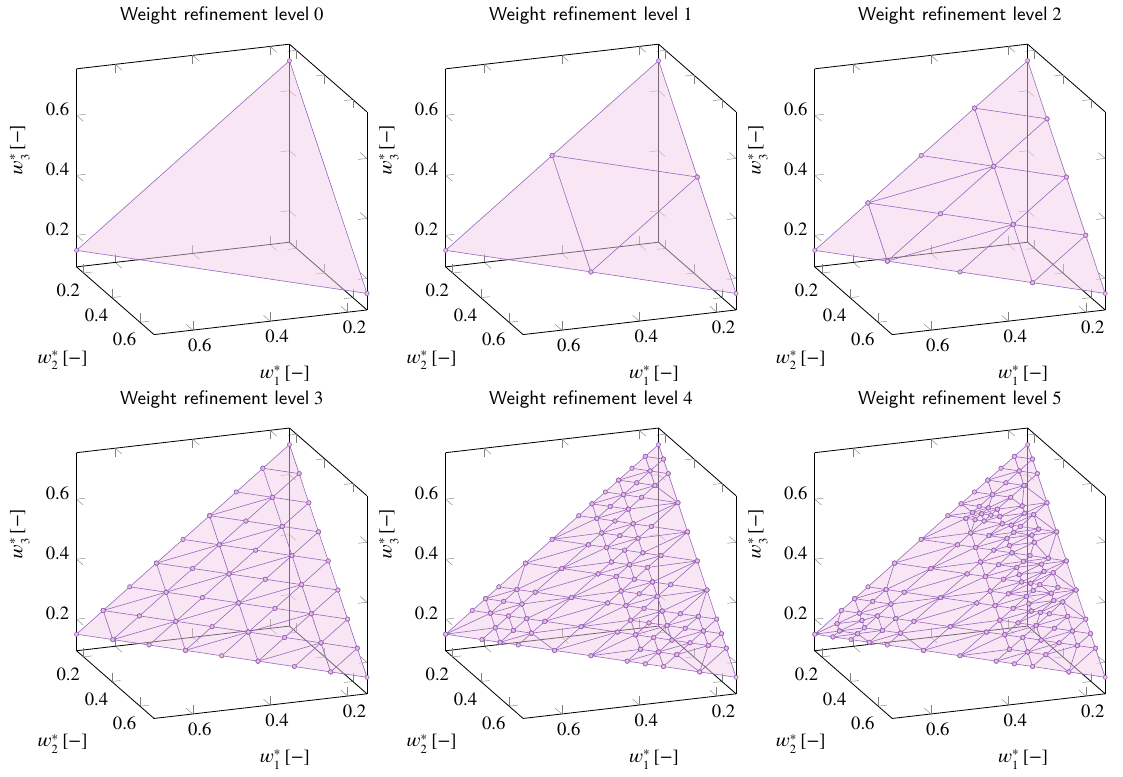}
    \caption{Weight refinement of the tri-objective minimum mean compliance problem.}
    \label{fig:CLBweightrefinement}
\end{figure}

\begin{figure}[ht]
    \centering
    \includegraphics[width=0.9\textwidth]{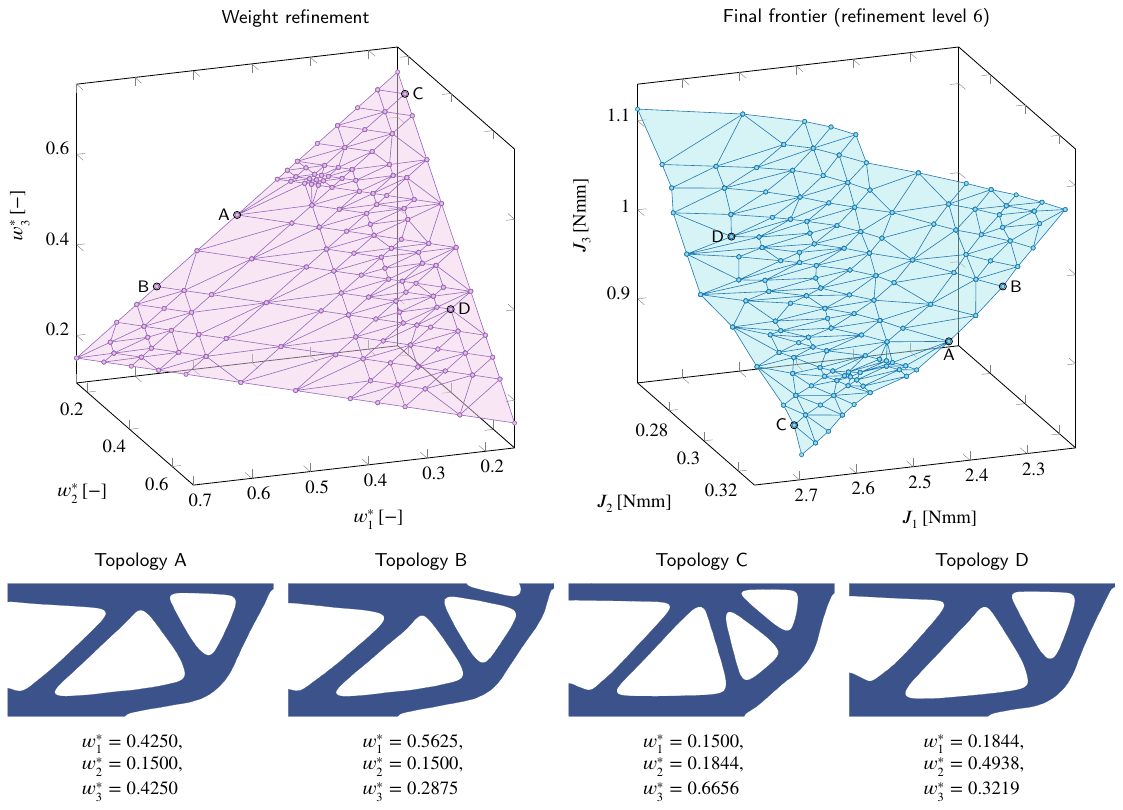}
    \caption{Final Pareto frontier approximation and selected topologies of the tri-objective minimum mean compliance problem.}
    \label{fig:CLBfrontierfinal}
\end{figure}

\begin{figure}[ht]
    \centering
    \includegraphics[width=0.9\textwidth]{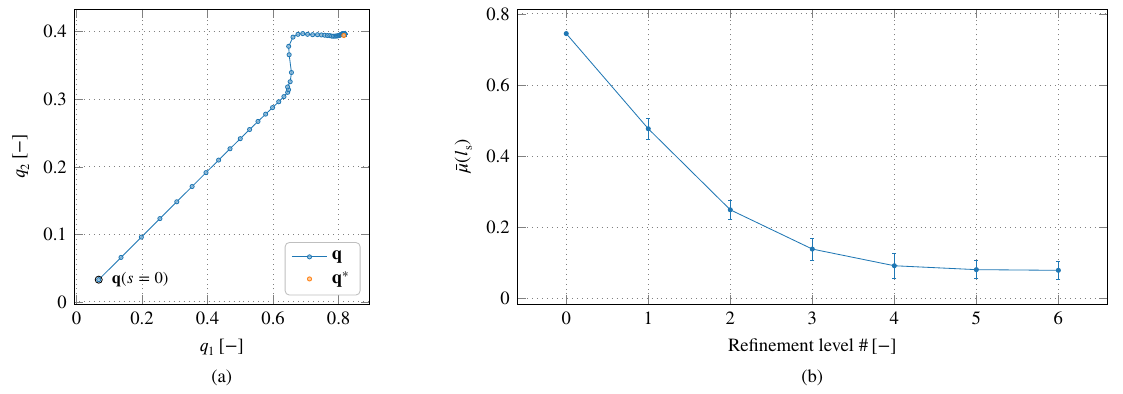}
    \caption{Evaluation plots of the tri-objective minimum mean compliance problem: (a) weight evolution of topology D, (b) convergence of the mean simplex length and its standard deviation.}
    \label{fig:CLBevaluation}
\end{figure}

As a final example, we consider a tri-objective minimum mean compliance problem ($\alpha \in \{1,2,3\}$) with varying loading and support conditions, see Fig.~\ref{fig:problemCLB}, to demonstrate the scalability of the proposed framework to higher-dimensional settings.
The corresponding optimization problem extends Problem~(\ref{eq:biminimummeancompliance}) by an additional objective functional and reads
\begin{customopti}|s|
    {vec\quad inf}{\phi\in H^2(\mathscr{D})}{
    \mathbf{J}[\mathbf{u}_1,\mathbf{u}_2,\mathbf{u}_3] = \left(J_1[\mathbf{u}_1],J_2[\mathbf{u}_2],J_3[\mathbf{u}_3]\right)^\top}{}{}{}
    \addConstraint{G[\Theta] = \frac{1}{V_0}\int_\mathrm{D}\Theta\,\mathrm{d}\Omega - V_\mathrm{f}}{\leq 0}{}{}
    \addConstraint{-\partial_j(\mathbb{C}_{ijkl}\varepsilon_{kl}(\mathbf{u}_\alpha))}{=0}{\quad\text{in}\quad\Omega}{}
    \addConstraint{u_{\alpha,i} n_i}{=0}{\quad\text{on}\quad\Gamma^\mathrm{r}_\alpha}{}
    \addConstraint{u_{\alpha,i}}{=0}{\quad\text{on}\quad\Gamma^\mathrm{u}}{}
    \addConstraint{(\mathbb{C}_{ijkl}\varepsilon_{kl}(\mathbf{u}_\alpha))n_j}{=t_i}{\quad\text{on}\quad\Gamma^\mathrm{t}_\alpha}{}
    \addConstraint{(\mathbb{C}_{ijkl}\varepsilon_{kl}(\mathbf{u}_\alpha))n_j}{=0}{\quad\text{on}\quad\partial\Omega\setminus(\Gamma^\mathrm{u}\cup\Gamma^\mathrm{r}_\alpha\cup\Gamma^\mathrm{t}_\alpha)}{}
    \label{eq:triminimummeancompliance}
\end{customopti}  
with objective functionals 
\begin{align*}
    J_\alpha[\mathbf{u}_\alpha] = \int_{\Gamma^\mathrm{t}_\alpha} t_i u_{\alpha,i} \,\mathrm{d}\Gamma.
\end{align*}
Since the derivation of the governing equations and sensitivity terms of Problem~(\ref{eq:triminimummeancompliance}) follows directly from the bi-objective case in Problem~(\ref{eq:biminimummeancompliance}) we focus on the numerical results in the following. The parameters are chosen as $C_{11} = C_{22} = 0.018$ and $B = 0.002$. The volume constraint is given as $V_\mathrm{f}=0.45$. For the weight evolution, we set $M_q = 0.5$, $B_q = 5$, and $K_q = 0.5$, while the maximum simplex length is prescribed as $l_{\mathrm{s,max}} = 0.15$. The initial weight discretization is given by
\begin{align*}
    \mathcal{W}_h^\ast=\{(0.70,0.15,0.15),(0.15,0.70,0.15),(0.15,0.15,0.70)\}.    
\end{align*}
Figure~\ref{fig:CLBfrontierrefinement} illustrates the iterative refinement of the Pareto frontier across successive refinement levels. 
Simplices whose maximum edge length exceeds the prescribed threshold $l_{\mathrm{s,max}}$ are refined accordingly. The associated refinement of the weight simplex is shown in Fig.~\ref{fig:CLBweightrefinement}.
A total of six refinement levels are performed to approximate the Pareto frontier. The final approximation and selected topologies are shown in Fig.~\ref{fig:CLBfrontierfinal}.
The convergence behavior in Fig.~\ref{fig:CLBevaluation}b indicates that the termination criterion is already achieved at refinement level $3$. The mean simplex edge length decreases rapidly during the first refinement steps and stabilizes thereafter, demonstrating convergence of the adaptive refinement procedure. At the same time, the standard deviation remains bounded, which suggests a uniform distribution of simplex sizes without the occurrence of degeneracies.
Figure~\ref{fig:CLBevaluation}a further shows the evolution of the weights for topology D from Fig.~\ref{fig:CLBfrontierfinal}, illustrating their adaptation to the local objective landscape.
  \section{Conclusions}
\label{sec:conclusions}
In this paper, we proposed a variational framework for multi-objective level set topology optimization. At the core of the approach is the interpretation of the level set function as a generalized coordinate of a fictitious material. By assigning suitable energy terms and applying Hamilton's principle, the evolution equation of the level set function can be derived in a multi-objective setting, resulting in a damped wave equation governing the design evolution. The objective functionals are combined through a weighted sum formulation and incorporated into a perturbation term that drives the system dynamics. An analysis of the derived equations naturally suggests a higher-level geometric interpretation of the problem. Under suitable regularity assumptions, the set of stationary solutions can be regarded as a geometric structure embedded in the space of objective functionals, within which the Pareto frontier is locally contained. In this interpretation, the weighting factors act as intrinsic coordinates. This viewpoint motivates the introduction of a dynamic evolution of the weights, allowing the optimization process to adapt to changes in the objective landscape during the evolution. As a result, a coupled dynamical system for the level set function and the weighting factors is obtained. The key findings and contributions of this work are summarized as follows:
\begin{enumerate}
    \item Development of a variational framework for multi-objective level set topology optimization, in which the level set function is interpreted as a generalized coordinate of a dynamical system. The corresponding evolution equation is derived from Hamilton's principle and results in a damped wave equation governing the design evolution.
    \item Introduction of a geometric interpretation of stationary solutions in the objective space. Under suitable regularity assumptions, the solution set can be regarded as a structured subset in which the Pareto frontier is locally embedded, with the weighting factors acting as intrinsic coordinates. This establishes a structure-aware perspective on multi-objective optimization.
    \item Development of an adaptive simplex decomposition strategy for the iterative approximation of the Pareto frontier, enabling targeted refinement in regions with higher curvature or stronger variation of the objectives.
    \item Introduction of a dynamic evolution of the weighting parameters, enabling an adaptive exploration of the objective landscape during the optimization process. The weighting dynamics are characterized by the parameters $M_q$, $B_q$, and $K_q$, which control responsiveness, damping, and scope of exploration. The conventional WS method arises as a limiting case of the proposed dynamics.
    \item Numerical results demonstrate that the proposed framework yields stable and uniform approximations of the Pareto frontier and remains effective for higher-dimensional objective spaces, as illustrated by several numerical examples, including a tri-objective optimization problem.
\end{enumerate}
The findings also open several avenues for future research. In particular, the explicit integration of the geometric structure of the solution set, together with enhanced simplex quality measures, directly into the evolution process represents a promising extension of the present framework. Furthermore, adaptive strategies for the automatic tuning of the weight evolution parameters are expected to further enhance the robustness and efficiency of the proposed approach.

  % ... 
  \section*{Acknowlegdement}
  This work was supported by the JSPS KAKENHI Grant Number 24KF0144.


  \bibliographystyle{elsarticle-num-names} 
  \bibliography{BIB_literature.bib}

\end{document}